\documentclass[3p]{elsarticle}

\usepackage{amsfonts,amsbsy,amssymb,amsmath}
\usepackage{lineno,hyperref}
\usepackage{graphicx}
\usepackage[table,xcdraw]{xcolor}
\usepackage{float}

\journal{Chemical Physics Letters}

\biboptions{compress}

\begin{document}

\begin{frontmatter}

\title{The Dynamical Significance of Valley-Ridge Inflection Points}

 \author[label1]{V. J. Garc\'ia-Garrido}\corref{mycorrespondingauthor}
 \ead{vjose.garcia@uah.es}
 \author[label2]{S. Wiggins}
 \ead{s.wiggins@bristol.ac.uk}

 \address[label1]{Departamento de F\'isica y Matem\'aticas, Universidad de Alcal\'a, \\ Madrid, 28871, Spain.\\[.2cm]}

 \address[label2]{School of Mathematics, University of Bristol, \\ Fry Building, Woodland Road, Bristol, BS8 1UG, United Kingdom.}

 \cortext[mycorrespondingauthor]{Corresponding author}

\begin{abstract}	
	In this paper we demonstrate that valley-ridge inflection (VRI) points of a potential energy surface (PES) have a dynamical influence on the fate of trajectories of the underlying Hamiltonian system. These points have attracted the attention of chemists in the past decades when studying selectivity problems in organic chemical reactions whose energy landscape exhibits a post-transition-state bifurcation in the region between two sequential saddles without an intervening energy minimum. To address the dynamical significance of valley-ridge inflection points, we construct a symmetric potential energy function that allows us to move the location of the VRI point while keeping the locations and energies of the critical points fixed. In this setup, we carry out a parametric study of the dynamical behavior of ensembles of trajectories in terms of the energy of the chemical system and the position of the VRI point. Our analysis reveals that the location of the VRI point controls the fraction of trajectories that recross the high energy saddle region of the PES without entering either of the potential wells that are separated by the low energy saddle.
\end{abstract}


\end{frontmatter}



\section{Introduction}
\label{sec:intro}

Chemical reactions exhibiting post-transition-state bifurcations (PTSBs) are a topic of current and growing interests in the organic chemistry community. In this context, after reaction, the system can evolve to two distinct products without passing through another transition state. An understanding to which product the system evolves (``selectivity'') offers the possibility of designing reactions with a desired outcome \cite{thomas2008,hong2014,Hornsby2014,hare2016}. 

An early review that catalogs a number of organic chemical reactions exhibiting PTSBs is \cite{ess2008}. More recent reviews \cite{rehbein2011,hare2017} describe a growing number of questions and directions for future investigations in this area. However, it is important to emphasize that a central theme for investigations in this topic is a need to understand the dynamics of organic reactions exhibiting PTSBs  \cite{tantillo2019wiggling}, and this, in turn, highlights the essential need for a phase space perspective of chemical reaction dynamics \cite{Agaoglou2019}. Phase space, which is the mathematical space comprised of the positions (the configuration space coordinates) and momenta of the underlying Hamiltonian system, provides the natural arena to explore dynamics. It is paramount to remark here that momentum is a crucial ingredient for the complete understanding of dynamics, and that, without it, there exists no way of explaining dynamical behavior of trajectories only from the topographical features of a PES. 

The basic features of two-dimensional potential energy surfaces describing this PTSB mechanism are well-known in the organic chemistry community, and are extensively described in the reviews cited above. Their topography typically displays four critical points: a high energy saddle, and a lower energy saddle separating two potential wells. In between the two saddle points there is a valley-ridge inflection point \cite{metiu1974,Valtazanos1986,Quapp1998,Quapp2004}, which is the point where the PES geometry changes from a valley to a ridge. The region between the two saddles forms a reaction channel and the dynamical issue of interest is how trajectories cross the high energy saddle, evolve towards the lower energy saddle, and select a particular well to enter. It is important to point out that the trajectories do not cross the lower energy saddle before entering one of the wells. Rather, the classical reaction path \cite{fukui1970} defined in terms of the landscape of the potential energy surface  bifurcates (i.e., splits into two curves) in a region shortly after the higher energy saddle \cite{quapp2004b}. This is consistent with the transition of this region of the PES from a valley to a ridge. and has encouraged the deeply rooted belief that VRI points play an important role in determining how trajectories choose which of the potential wells to enter. However, detailed trajectory studies on potential energy surfaces of this type, often referred to as ``VRI potential energy surfaces'', have not revealed a definitive role for the VRI point in the dynamical evolution of trajectories \cite{collins2013,katsanikas2020PRE,gg2020cplett,makrina2020cplett}, even in the context of quantum wavepacket dynamics \cite{lasorne2003,lasorne2005}. The goal of this paper is to demonstrate for the first time that VRI points have a dynamical effect on trajectories, and our results give a positive answer to this question. This is an important and rather surprising result, since VRI points typically are not equilibrium points of Hamilton's equations of motion.

Reacting trajectories crossing the high energy saddle have three possible fates when studied at short to moderate timescales. They can enter one well or the other and we refer to these as the ``top'' and ``bottom'' wells. This terminology will be made clear when we explicitly define and sketch the potential energy surface in Section \ref{sec:sec1}. Another possibility is that they can return to where they came from and ``recross'' the region of the higher energy saddle without entering either well. The branching ratio is the ratio of the total number of trajectories that enter the top or bottom wells. This number quantifies the notion of selectivity. For our study, we will use a symmetric VRI potential energy surface, i.e., one where the top and bottom wells are symmetrically related. In this case, equal numbers of trajectories enter the top and bottom wells. Moreover, our PES has the property that the location of the VRI point can be moved along the line that connects the high energy and low energy saddles without affecting the symmetry of the PES and the locations and energies of its critical points. This setup has the advantage of allowing us to probe the dynamical effect of the location of the VRI point on trajectories in a systematic way where the branching ratio is not affected. In this setting, we show that for symmetric VRI potentials the location of the VRI point directly affects the number of recrossing trajectories. 

This paper is outlined as follows. In Section \ref{sec:sec1} we introduce the PES model that we have used to address the dynamical influence of VRI points on the fate of trajectories for the underlying Hamiltonian system. We also describe the experimental setup developed to test numerically this effect on ensembles of trajectories that initially move across the high energy saddle. Next, Sec. \ref{sec:sec2} is devoted to describing the results of this work. We show by means of running ensembles of trajectories how the location of the VRI point controls the fraction of trajectories that, after crossing the high energy saddle, they recross it without entering either of the potential wells. We do so by calculating fate maps and performing an statistical analysis from the number of trajectories that display distinct dynamical behaviors. Finally, in Sec. \ref{sec:conc} we present the conclusions of this work and briefly discuss some questions that we will pursue in the near future to extend this research further.


\section{The Potential Energy Surface Model and Experimental Setup}
\label{sec:sec1}

In this section we describe the PES that we have devised to study the dynamical impact of VRI points on trajectories. Much effort has been devoted in the chemistry community during the past years to the task of designing PES with symmetric and asymmetric PTSB regions, and also to the analysis of how such topographical features affect selectivity in chemical reactions \cite{chuang2020}. However, trying to understand dynamics in complex situations where many factors of the PES are varied simultaneously, such as the location and energies of the critical points, can make this endeavour a challenging task that obscures the fundamental underlying mechanisms at play. We have decided to follow in this paper a bottom-up strategy to address the question of how VRI points have dynamical significance. To do so, we construct a simplified PES model inspired in the work carried out in \cite{collins2013}. The advantage that our PES model brings is two-fold. First, since we will work with a symmetric PES, we know that the branching ratio of trajectories always remains equal. But the most important characteristic of our model system is that we can move the location of the VRI point along the line that connects both saddles without affecting the energies and locations of the remaining critical points of the PES. This allows us to perform a systematic analysis of how the location of the VRI point affects the evolution of ensembles of trajectories, and we study this question in terms of the energy of the system. 

\begin{figure}[htbp]
	\begin{center}
		A)\includegraphics[scale = 0.24]{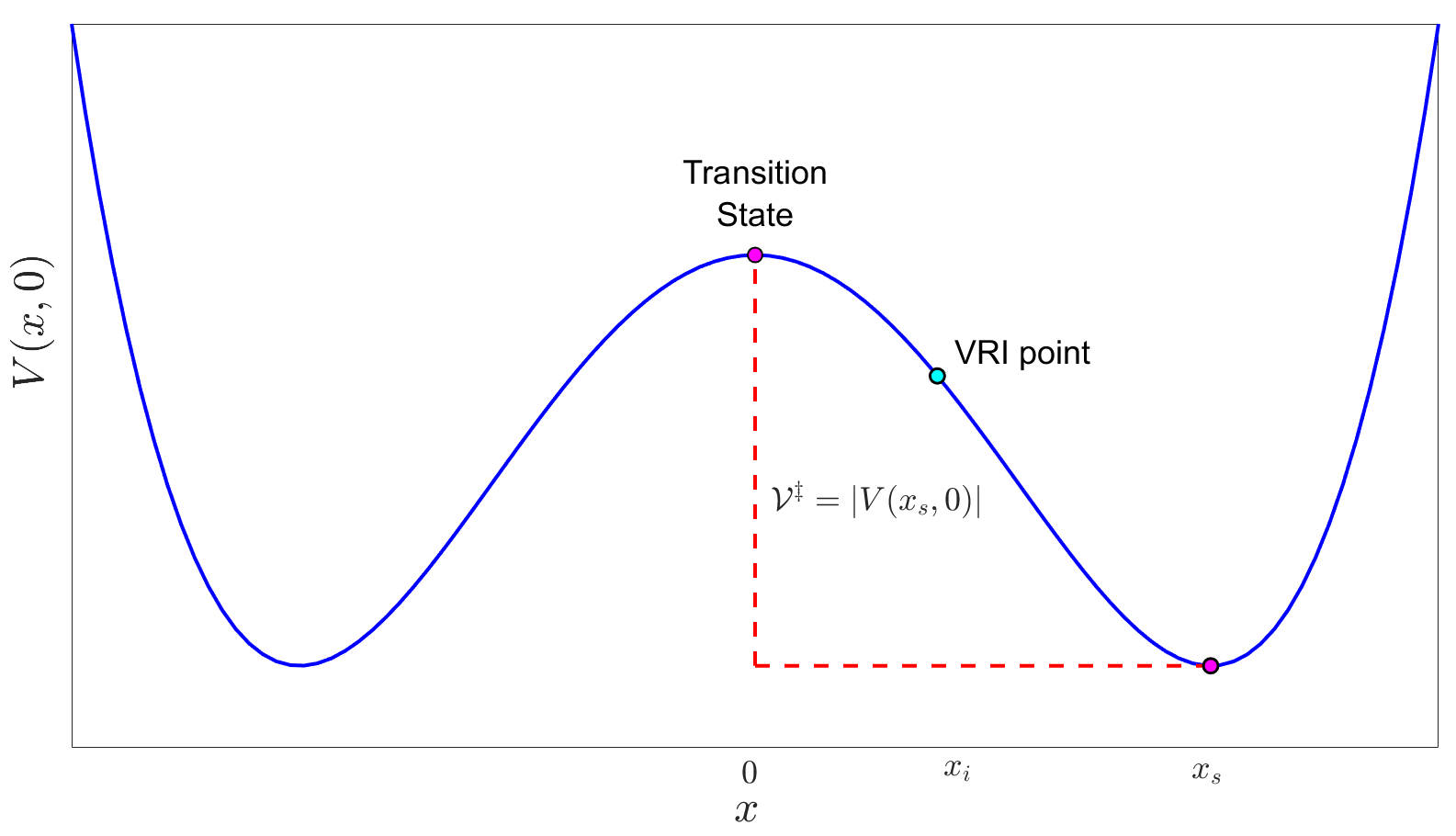} \\[.2cm]	
		B)\includegraphics[scale = 0.27]{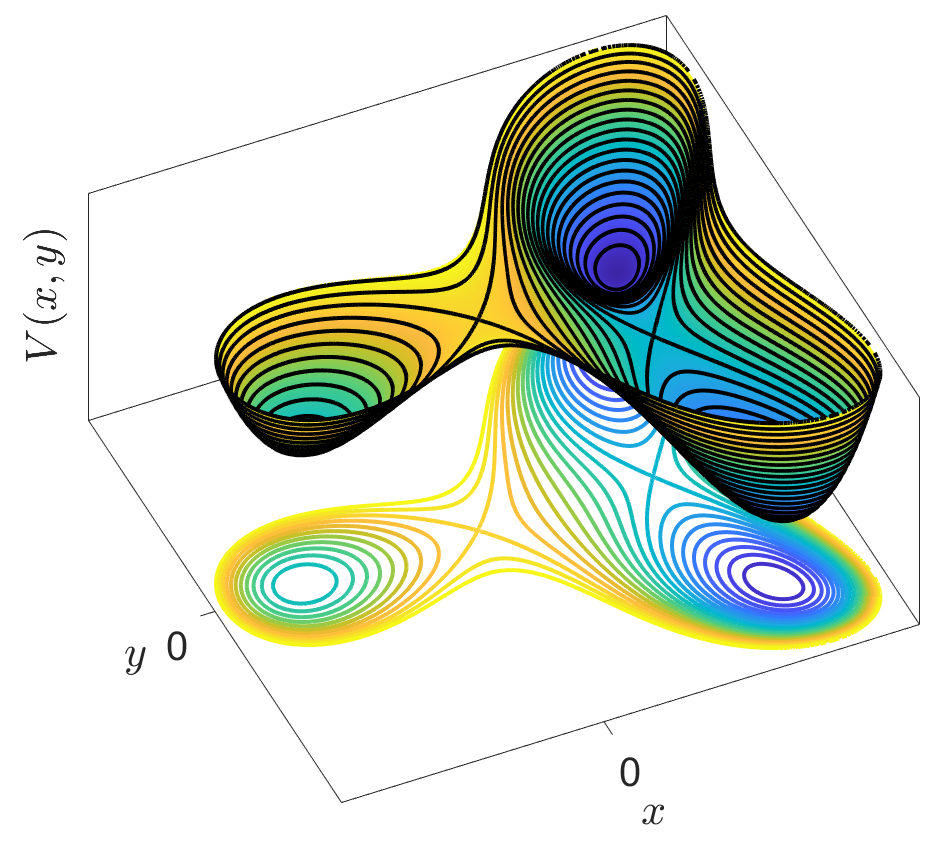}
		C)\includegraphics[scale = 0.28]{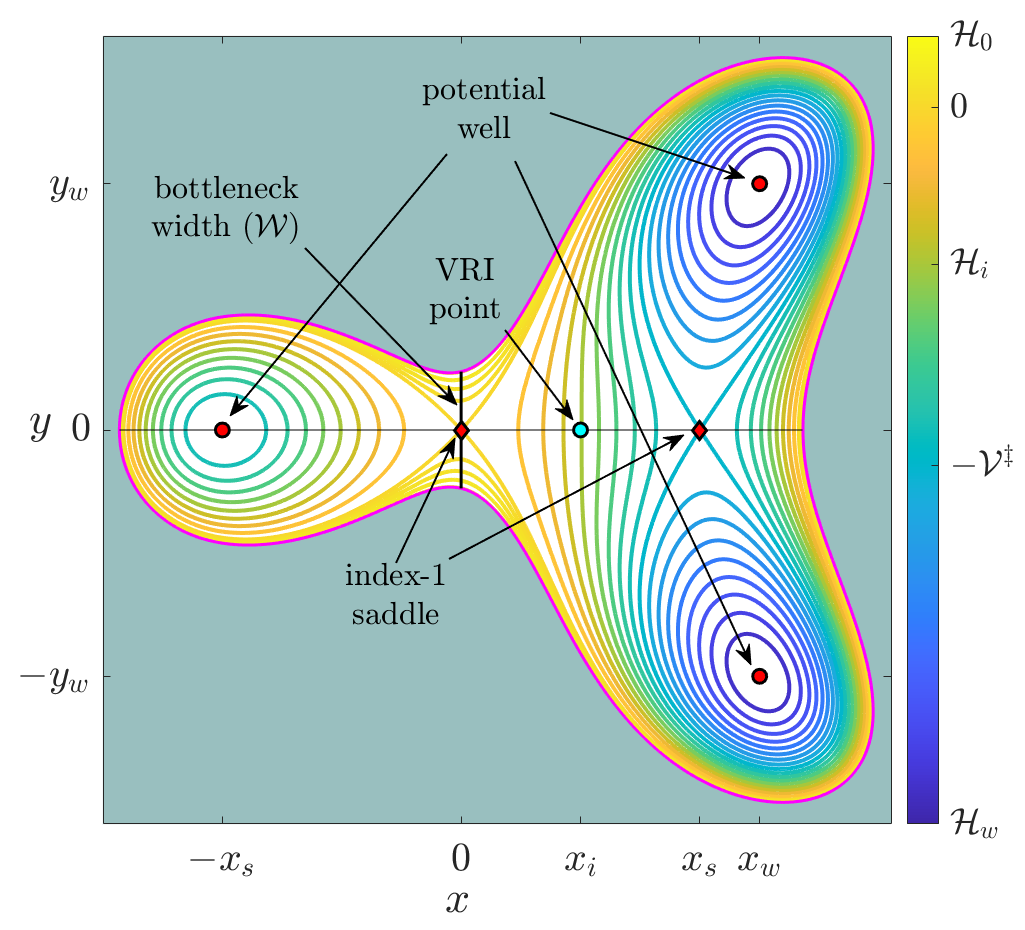}
	\end{center}
	\caption{A) Section of the PES described in Eq. \eqref{eq:pes} along the $x$-axis. We have marked the location of the saddles as magenta dots. B) Three-dimensional representation of the potential energy surface. C) Equipotential contours, locations and energies of the critical points and of the VRI point of the potential energy function.}
	\label{fig:pes}
\end{figure}

We construct an energy landscape consisting of three potential wells and two index-1 saddle points. One of the saddle points, which is located at the origin and is the critical point with the highest energy in the system, separates the PES into two regions. On the left, we have one well that corresponds to reactants, while on the right there are two product wells separated by a lower energy saddle. In this setting, and in order to simplify the analysis further, we will impose that the PES is symmetric with respect to the $x$-axis and also that the two saddles of the system lie on the $x$-axis. This condition implies in particular that there exists a VRI point between them. A representation of the topography of this model PES is displayed in Fig. \ref{fig:pes}. Recall that in the vicinity of a VRI point, the intrinsic reaction coordinate bifurcates due to the shape of the PES \cite{quapp2004b,birney2010} and this gives rise to a reaction mechanism known as a two-step-no-intermediate mechanism \cite{singleton2003}. Mathematically, at a VRI point two conditions are met: the Gaussian curvature of the PES is zero, which implies that the Hessian matrix has a zero eigenvalue, and also the gradient of the potential is perpendicular to the eigenvector corresponding to the zero eigenvalue. Geometrically, this means that the landscape of the PES in the neighborhood of the VRI changes its shape from a valley to a ridge. Mathematically, above conditions can be written as:
\begin{equation}
	\begin{cases}
		\det\left(\text{Hess}_V\right) = 0 \\[.2cm]
		\left(\nabla V\right)^T \text{adj}\left(\text{Hess}_V\right) \, \nabla V = 0
	\end{cases} \,,
	\label{eq:vri_conds}
\end{equation}
where $V$ is the potential energy function, the term $\text{Hess}_V$ corresponds to the Hessian matrix and $\text{adj}\left(\text{Hess}_V\right)$ represents the adjugate matrix of the Hessian of the PES. It is important to remark here that VRI points are not critical points of the PES. Despite this fact, we show in this paper that they play a relevant dynamical role that has a measurable and quantifiable influence on the evolution and fate of the system trajectories.

Consider a Hamiltonian system with two degrees-of-freedom (DoF) defined as the classical sum of kinetic and potential energy in the form:
\begin{equation}
	\mathcal{H}(x,y,p_x,p_y) = \dfrac{p_x^2}{2 m_1} + \dfrac{p_y^2}{2 m_2} + V(x,y) \;,
\end{equation}
where $m_1$ and $m_2$ are the masses associated to the $x$ and $y$ DoF, respectively, and the potential energy function is given by the expression:
\begin{equation}
	V(x,y) = \dfrac{\mathcal{V}^{\ddagger}}{x_s^4} x^2 (x^2 - 2x_s^2) + A y^2 (x_i - x) + y^4 (B + Cx)
	\label{eq:pes}
\end{equation}
where $A$, $B$, $C$ are free model parameters, $x_s$ is the $x$-coordinate location of the index-1 saddle that sits between the two potential wells on the right hand side of the origin, $\mathcal{V}^{\ddagger}$ represents the energy barrier height of the index-1 saddle at the origin, and $x_i$ denotes the $x$-coordinate of the VRI point that exists at the configuration space point $(x_i,0)$ between the high energy saddle at the origin and the low energy saddle located at $(x_s,0)$. We depict the geometry of the PES given by Eq. \eqref{eq:pes} in Fig. \ref{fig:pes}. In panel A) we show the profile of the potential energy function along the $x$-axis, and plots B) and C) display the map of equipotential contours together with the locations and energies of the critical points, and also of the VRI point. It is important to note here that the model parameters $A$, $B$ and $C$ are chosen by means of solving a linear system of equations that results from imposing that the two potential wells on the right hand side of the origin are located at the points $(x_w,\pm y_w)$ and have an energy $\mathcal{H}_w$.

Hamilton's equations of motion that govern the dynamics of trajectories for this system are given by:
\begin{equation}
	\begin{cases}
		\dot{x} = \dfrac{\partial \mathcal{H}}{\partial p_x} = \dfrac{p_x}{m_1} \\[.5cm]
		\dot{y} = \dfrac{\partial \mathcal{H}}{\partial p_y} = \dfrac{p_y}{m_2} \\[.5cm]
		\dot{p}_x = -\dfrac{\partial V}{\partial x} = \dfrac{4\mathcal{V}^{\ddagger}}{x_s^4}x\left(x_s^2 - x^2\right) + Ay^2 - Cy^4  \\[.5cm]
		\dot{p}_y = -\dfrac{\partial V}{\partial y} = 2Ay\left(x-x_i\right)-4y^3\left(B+Cx\right)
	\end{cases}
	\label{eq:hameq}
\end{equation}
and for our analysis we will use the values $m_1 = m_2 = 1$ for the masses of the DoF. For this 2 DoF Hamiltonian, dynamics takes place in a four-dimensional phase space and, since energy is conserved, motion is constrained to a three-dimensional energy surface. It is a simple exercise to show that the eigenvalues of the Jacobian matrix evaluated at the equilibrium point at the origin, which characterize the linearized dynamics in its neighborhood, are given by:
\begin{equation}
	\lambda_{1,2} = \pm \dfrac{2 \sqrt{\mathcal{V}^{\ddagger}}}{x_s} \quad,\quad \lambda_{3,4} = \pm \sqrt{2A x_i} \, i \;.
\end{equation}
From these expressions we can clearly see that the location of the VRI point, described by the variable $x_i$, has a direct effect on the linearized angular frequency of vibration in the bottleneck region of the index-1 saddle point at the origin. We can estimate the configuration space width of the bottleneck region about the transition structure at the origin by the following procedure. Given an energy level for the system $\mathcal{H} = \mathcal{H}_0 > 0$ above that of the high energy saddle, take the vertical line in configuration space that passes through the origin and connects the corresponding two equipotential curves with energy $V =  \mathcal{H}_0$, see Fig. \ref{fig:pes} C). The bottleneck width is given by following formula:
\begin{equation}
	\mathcal{W} = 2 \sqrt{-\dfrac{A x_i}{2B} + \sqrt{\left(\dfrac{A x_i}{2B}\right)^2 + \dfrac{\mathcal{H}_0}{B}}} \;,
\end{equation} 
showing that the bottleneck width is also a function of the VRI location.

We finish this section by describing the setup for the numerical experiments carried out  for this model system. For our simulations we will use the following model parameters. The barrier height at the origin is set to $\mathcal{V}^{\ddagger} = 0.5$, the lower energy saddle is at $x_s = 1$, and the potential wells on the right side of the origin are located at the coordinates $(x_w,y_w) = (1.25,\pm 1)$ with energy $\mathcal{H}_w = -1$. In order to probe the dynamical influence of the VRI point we have run ensembles of trajectories taken from two different configurations. First, we have uniformly sampled the vertical line in configuration space that passes through the saddle point at the origin and connects the two equipotentials with the same total energy of the system $\mathcal{H} = \mathcal{H}_0$. We define the initial conditions by extending these configuration space points to the full phase space, and to do so, we initialize the trajectories  with all the momentum allocated along the $x$-coordinate in the positive direction. The physical interpretation of this condition is that all the trajectories initially cross the high energy saddle horizontally from left to right, that is, we set $p_y = 0$. Moreover, the $p_x$ component of momentum has to be selected so that the initial condition satisfies the energy constraint. This set of points can be written as:
\begin{equation}
	\mathcal{C}(\mathcal{H}_0) = \left\lbrace (x,y,p_x,p_y) \in \mathbb{R}^4 \,:\, x = p_y = 0 \;,\; V(0,y) \leq \mathcal{H}_0 \;,\; p_x = \sqrt{2\left[\mathcal{H}_0 - V(0,y) \right]} \right\rbrace
	\label{eq:cs_slice}
\end{equation}
The other set that we will use to sample initial conditions is defined by the phase space slice that orthogonally intersects the $x$-axis at the origin, and we consider that $p_x > 0$. This assumption represents physically the situation where trajectories initially evolve by entering the PES region that contains the two symmetrically related potential wells to the right of the high energy saddle. Notice that in this case we allow that some of the initial momentum of the trajectory can be directed along the $y$ coordinate, that is $p_y \neq 0$. Mathematically, we can write this phase space plane of initial conditions as:
\begin{equation}
	\mathcal{P}(\mathcal{H}_0) = \left\lbrace (x,y,p_x,p_y) \in \mathbb{R}^4 \,:\, x = 0 \;,\; p_x(y,p_y;\mathcal{H}_0) > 0 \right\rbrace
	\label{eq:ps_slice}
\end{equation}

\section{Results}
\label{sec:sec2}

\begin{figure}[htbp]
	\begin{center}
		A)\includegraphics[scale = 0.22]{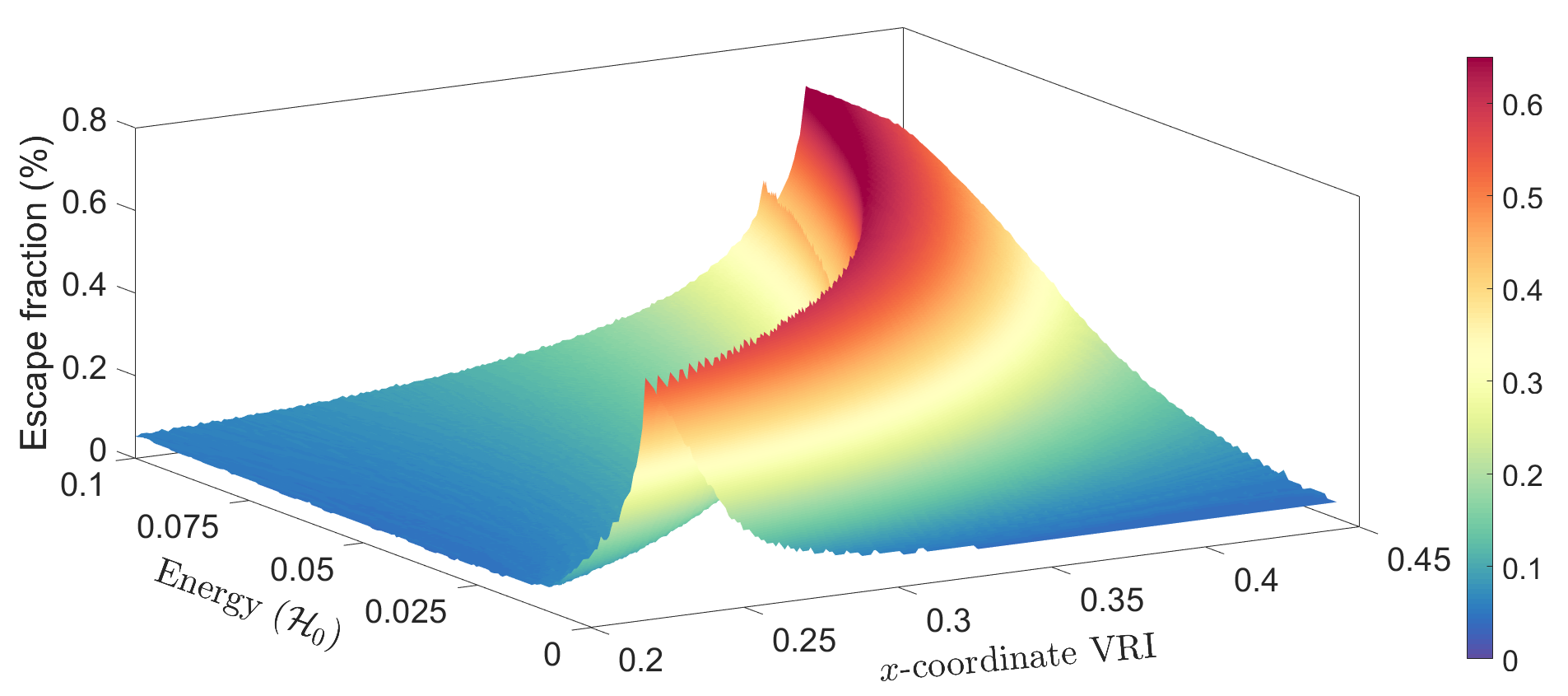}
		B)\includegraphics[scale = 0.22]{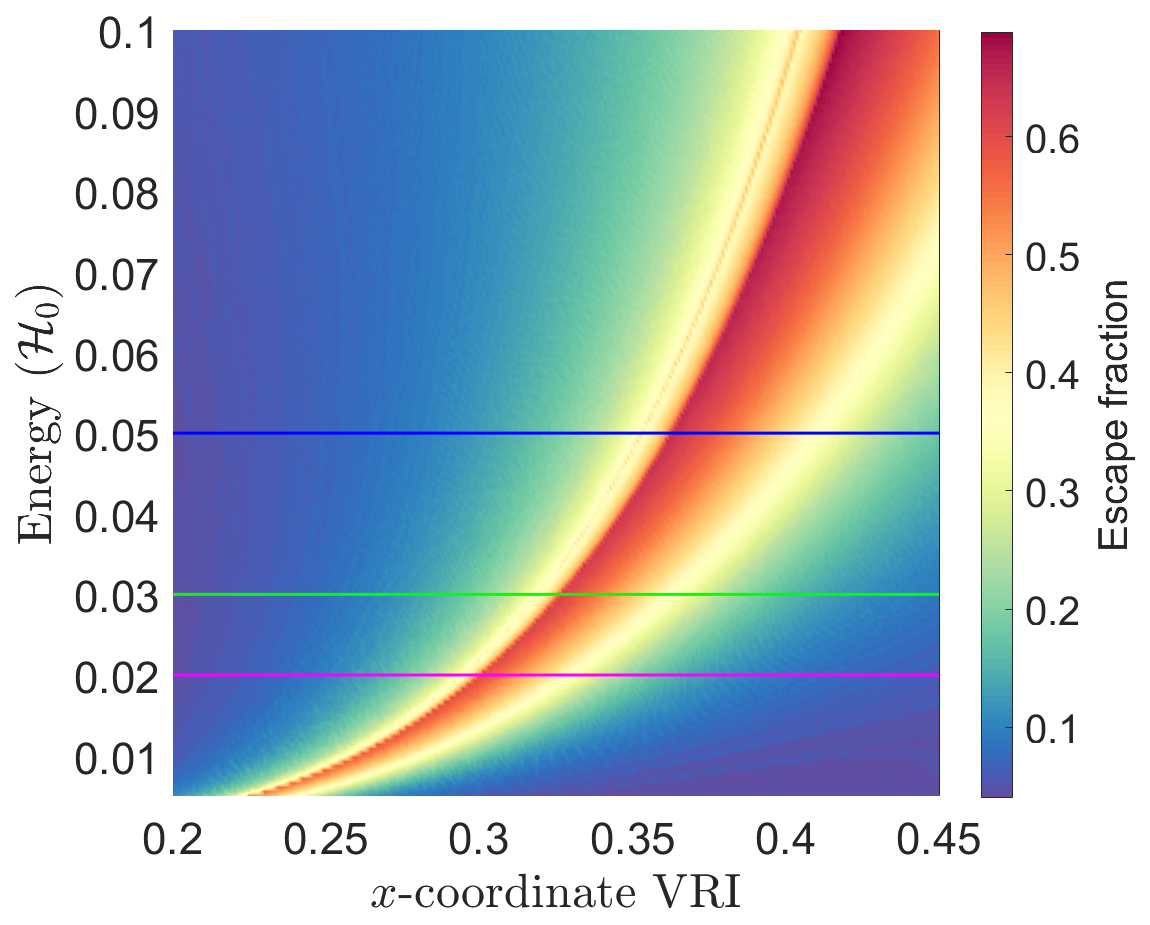}
		C)\includegraphics[scale = 0.22]{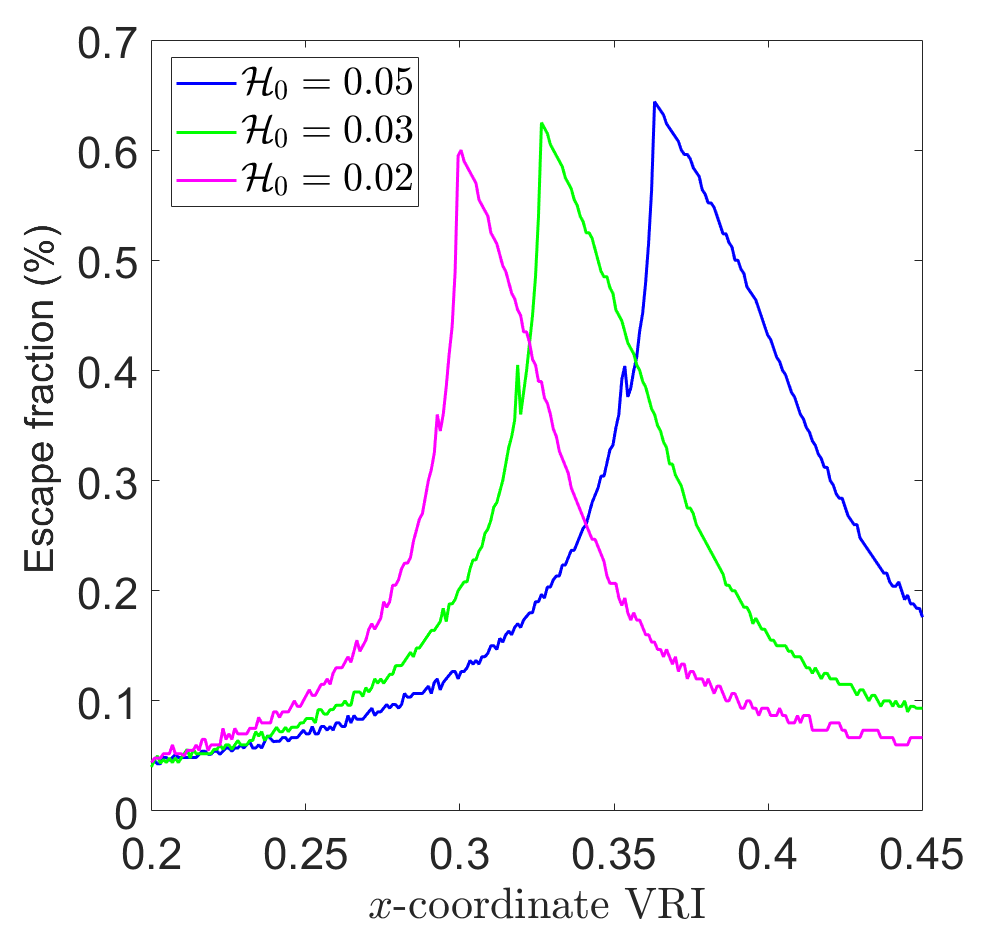}
	\end{center}
	\caption{A) Fraction of recrossing (escaping) trajectories as a function of the systems' energy and the location of the VRI point. Statistics are calculated for a uniform ensemble of initial conditions selected on the configuration space line described in Eq. \eqref{eq:cs_slice}. B) Top view of panel A). C) Recrossing fraction for different energy levels marked in B).}
	\label{vri_energy}
\end{figure}

We begin our analysis of the impact that VRI points have on the dynamics of trajectories by launching an ensemble of initial conditions uniformly sampled along the configuration space line $\mathcal{C}(\mathcal{H}_0)$ in Eq. \eqref{eq:cs_slice} located at the bottleneck region of the high energy saddle. We will use a linear density of 500 trajectories per unit length of bottleneck width $\mathcal{W}$, and this numerical experiment is carried out for a range of energies of the system from $\mathcal{H}_0 = 0.005$ to $\mathcal{H}_0 = 0.1$ with a step of $\Delta \mathcal{H}_0 = 0.001$. We also look at a range of values for the location of the VRI point from $x_i = 0.2$ to $x_i = 0.45$ with a step of $\Delta x_i = 0.0025$. In order to classify those trajectories that enter the top or bottom well regions, or those that recross the high energy saddle (they escape without giving rise to products), we impose the following condition. Those trajectories that get close to the top well (resp. bottom well) and enter a circle of radius $R = 0.2$ centered at the top well (resp. bottom well) critical point are stopped and counted accordingly. The results of this parametric study are shown in Fig. \ref{vri_energy}, where we display the fraction of recrossing trajectories as a function of energy and of the VRI point location. This analysis reveals that for every energy level of the system, there is a maximum peak of recrossing trajectories for a given value of the VRI point location. Interestingly, the recrossing fraction obtained can attain very large values, more than $60\%$ of the trajectories recross the saddle region at the origin, even for small energies. Moreover, to the right and left of the maximum value, the recrossing fraction decreases almost monotonically, although some fluctuations in its value are clearly observed. In fact, we can also identify in the plots two distinctive sharp peaks, one where the recrossing fraction is maximal, and another located to the left of this ridge. It is important to remark here that similar dynamical behavior where a significant percentage of recrossing trajectories has been observed and reported in chemical reactions such as in the Diels-Alder reaction, see \cite{singleton2009}.

\begin{figure}[htbp]
	\begin{center}
		A)\includegraphics[scale = 0.21]{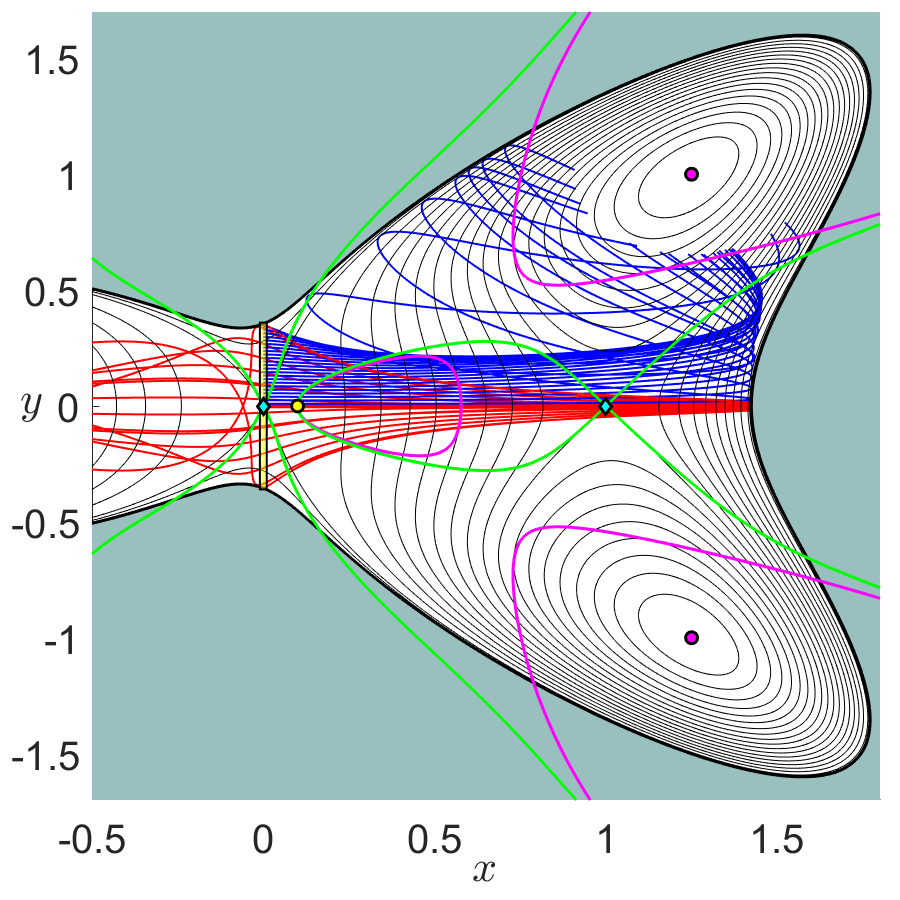}
		B)\includegraphics[scale = 0.21]{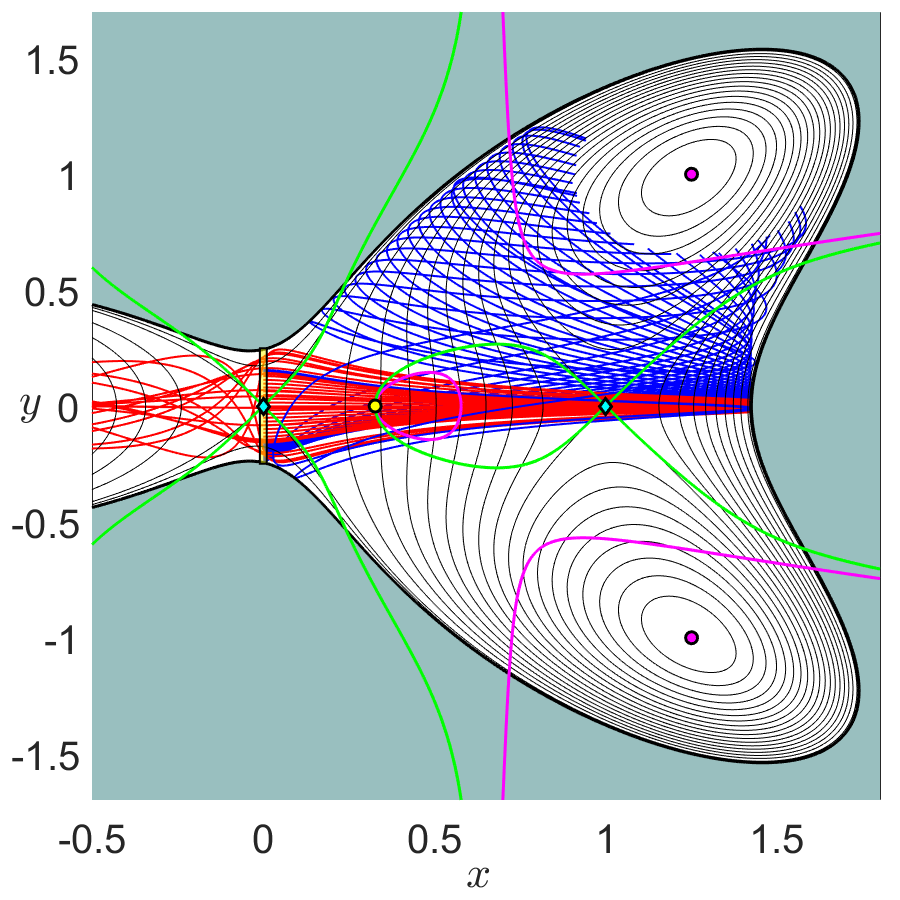}
		C)\includegraphics[scale = 0.21]{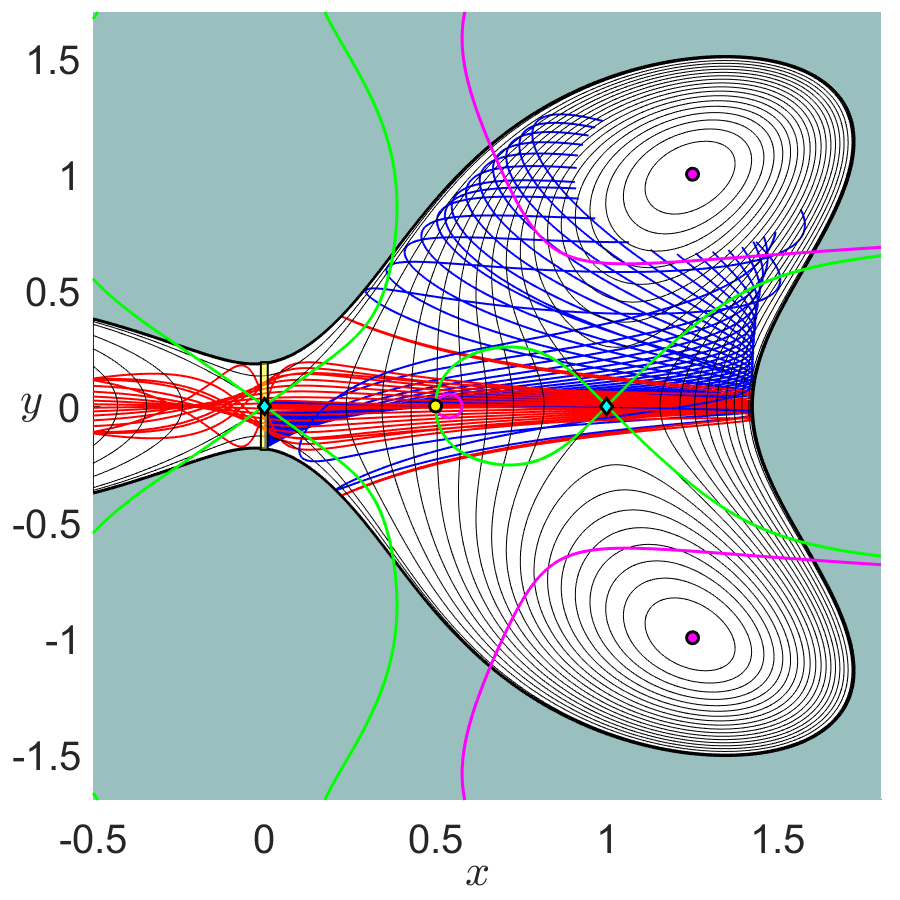}
	\end{center}
	\caption{Trajectory evolution of an ensemble of initial conditions selected on the configuration space set given by Eq. \eqref{eq:cs_slice} (yellow rectangle) for an energy $\mathcal{H}_0 = 0.03$. The location of the VRI point is different in all the panels: A) The VRI point is at $x_i = 0.1$. B) $x_i = 0.3265$. C) $x_i = 0.5$.		 Recrossing trajectories are depicted in red while those entering the top well region of the PES are displayed in blue. For clarity of the plots, we have omitted those trajectories that go to the bottom well because of the symmetry of the PES. The location of the saddles, the VRI point and the potential wells are also marked with cyan diamonds, yellow and magenta circles respectively. We have also overlaid the curves that represent the two conditions that VRI points satisfy as described in Eq. \eqref{eq:vri_conds}. The magenta curve corresponds to values where the Hessian matrix has zero determinant, while the green curve depicts the condition for the adjugate of the Hessian matrix.}
	\label{configSp_evol}
\end{figure}

We take a look next at how the trajectories behave and, in order to get an idea of the regions of the PES they traverse along their evolution, we depict their projection onto configuration space. Consider an energy level of $\mathcal{H}_0 = 0.03$, we would like to compare how the ensemble evolves for three different values of the VRI point location. For this purpose we select the value where recrossing is maximum, and this occurs at $x_i = 0.3265$, and two other values $x_i = 0.1$ and $x_i = 0.5$ on either side of the maximum peak for which the recrossing fraction decreases substantially. Along the line $\mathcal{C}(\mathcal{H}_0)$ in Eq. \eqref{eq:cs_slice} we select a uniform ensemble of initial conditions and simulate the trajectories until they enter the region of either of the wells or until they recross the high energy saddle, whatever happens first. We plot in Fig. \ref{configSp_evol} the results of this analysis. It is interesting to point out here that the VRI point seems to be having a lensing effect on the trajectory ensemble, focusing the trajectories on the wall of the PES opposite to the saddle at the origin. As the VRI point reaches the location $x_i = 0.3265$, this focusing mechanism is enhanced and becomes 'maximal', and most of the trajectories that bounce off the wall recross the high energy saddle region, going back to where they started, entering the reactant well on the left of the PES. Notice also the caustic-like pattern formed by the recrossing trajectories on the of the PES. This type of trajectory behavior has been reported in other studies concerned with how trajectories moving in a force field escape from a given region \cite{maier1993}.

\begin{figure}[htbp]
	\begin{center}
		A)\includegraphics[scale = 0.36]{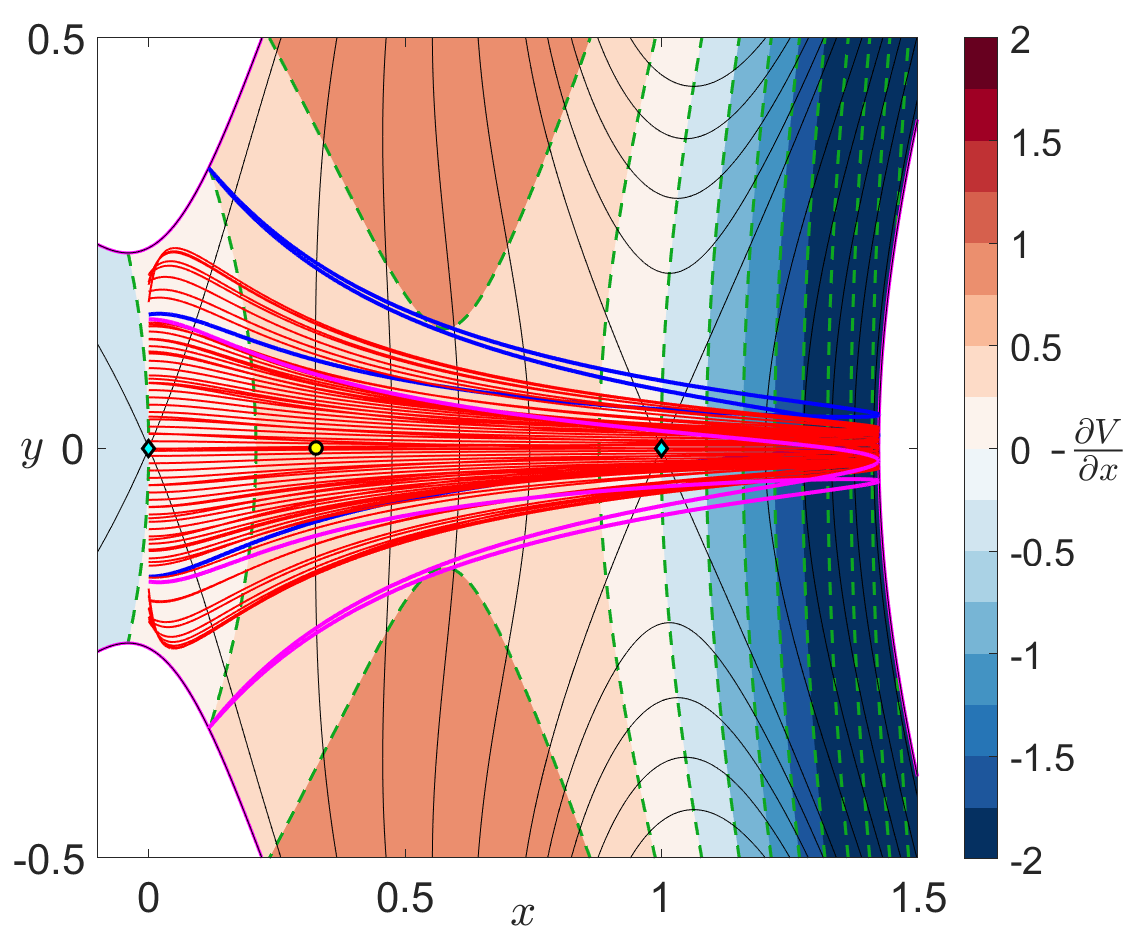}
		B)\includegraphics[scale = 0.36]{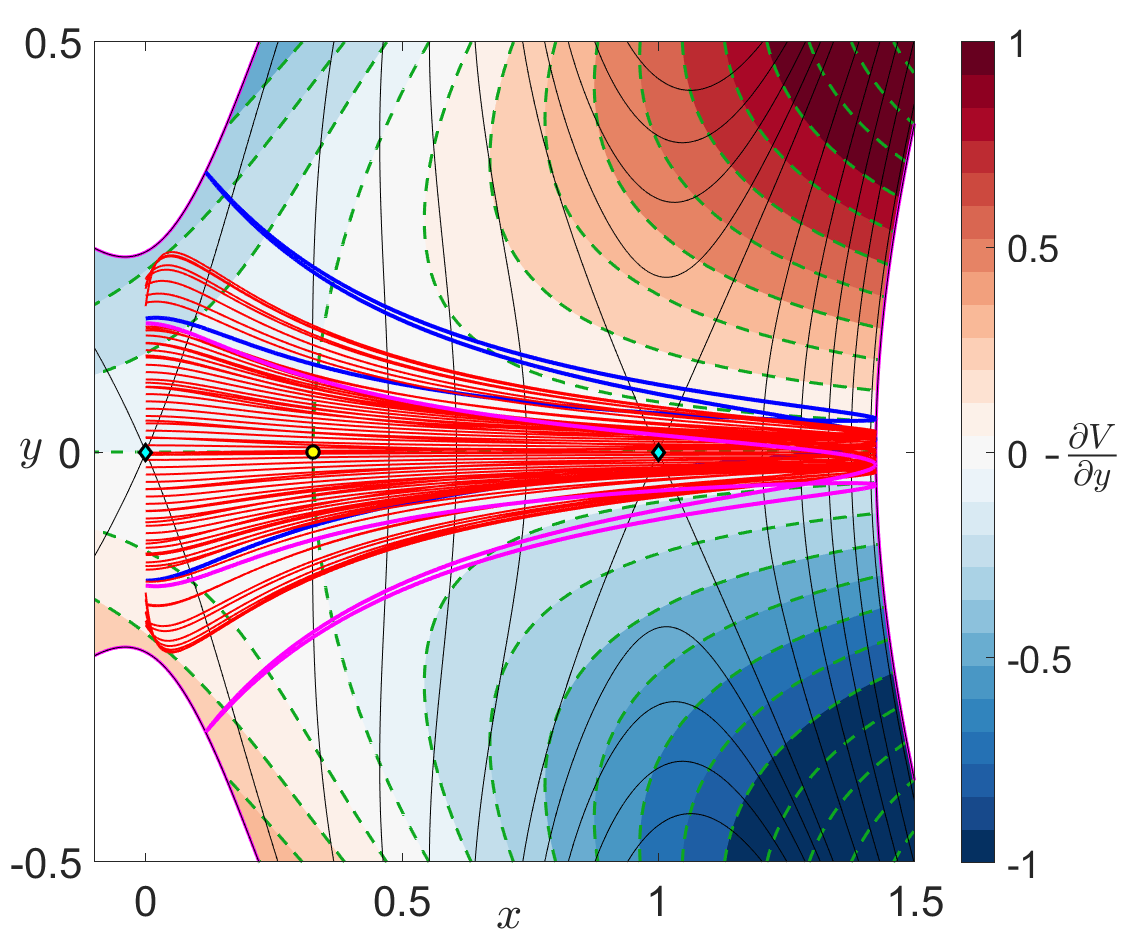}
	\end{center}
	\caption{Evolution of recrossing trajectories, initialized on the line $\mathcal{C}\left(\mathcal{H}_0\right)$ given in Eq. \eqref{eq:cs_slice}, projected onto configuration space. The energy of the system is $H_0 = 0.03$ and the VRI point is located at $x_i = 0.3265$. The saddles and the VRI point are marked as cyan diamonds and a yellow circle respectively. In panels A) and B) we have overlaid the plot with the contour levels of the components of the force at each point of the PES. The blue and magenta curves depict limiting trajectories, that is, they represent those recrossing trajectories that get 'closest' to entering the top and bottom well regions of the PES respectively.}
	\label{configSp_evol2}
\end{figure}

\begin{figure}[htbp]
	\begin{center}
		A)\includegraphics[scale = 0.2]{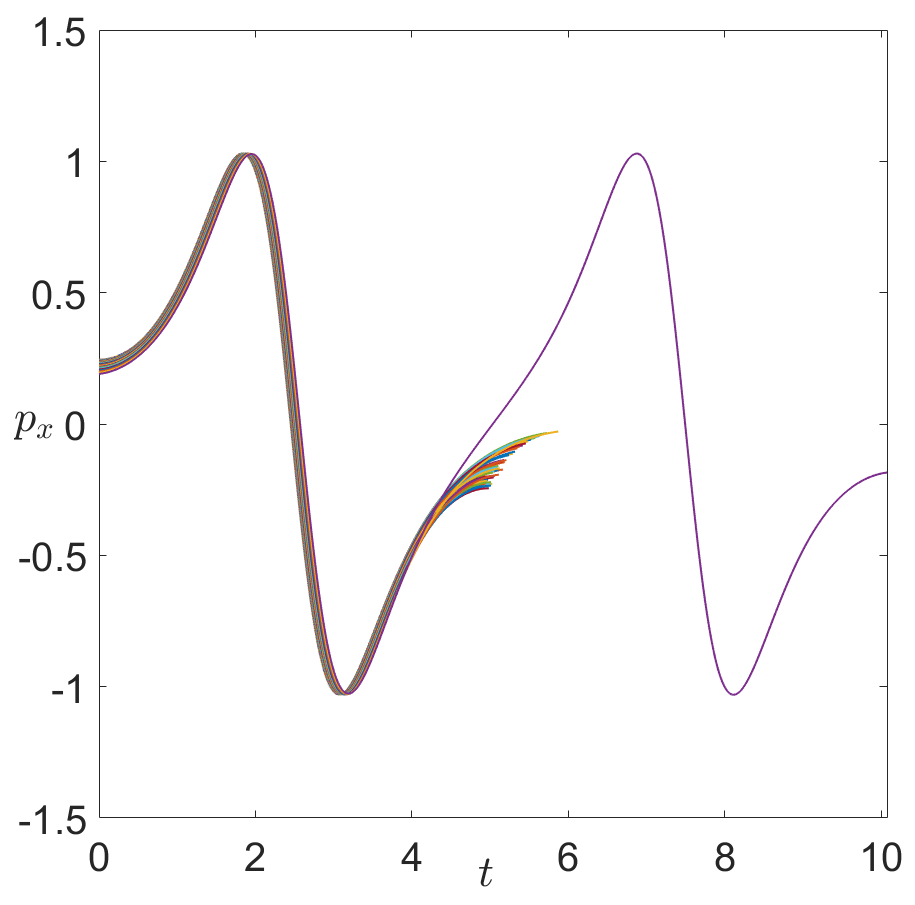}
		B)\includegraphics[scale = 0.2]{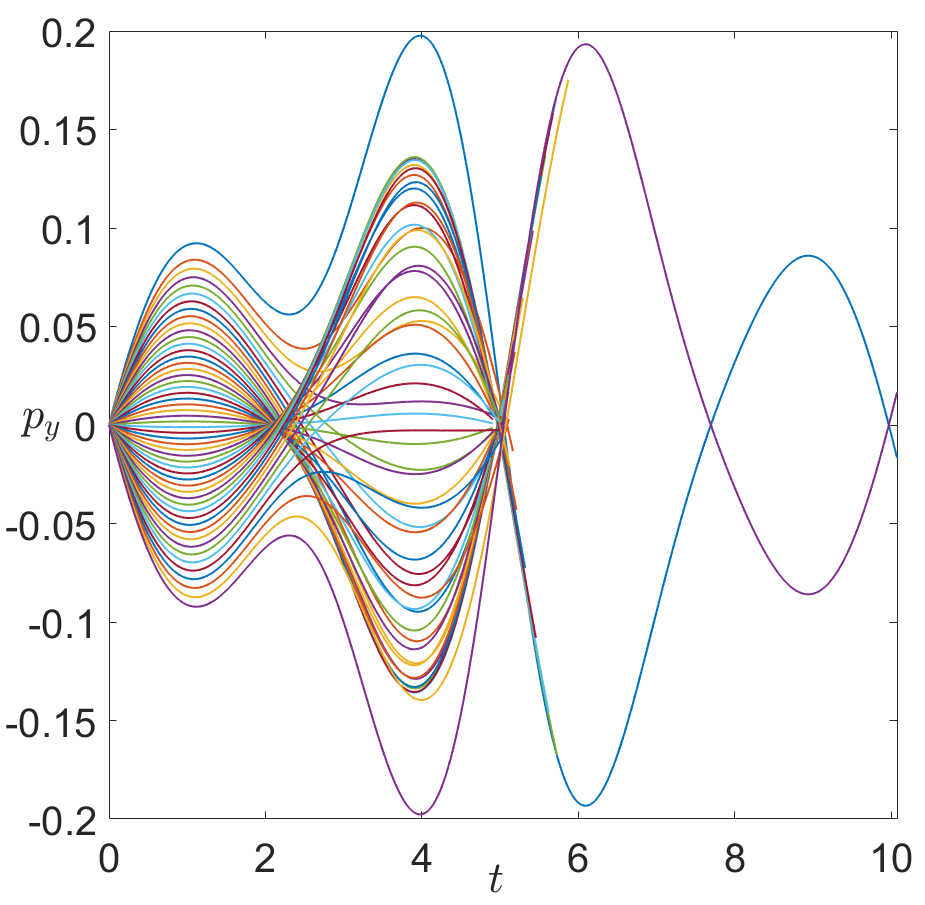}
		C)\includegraphics[scale = 0.2]{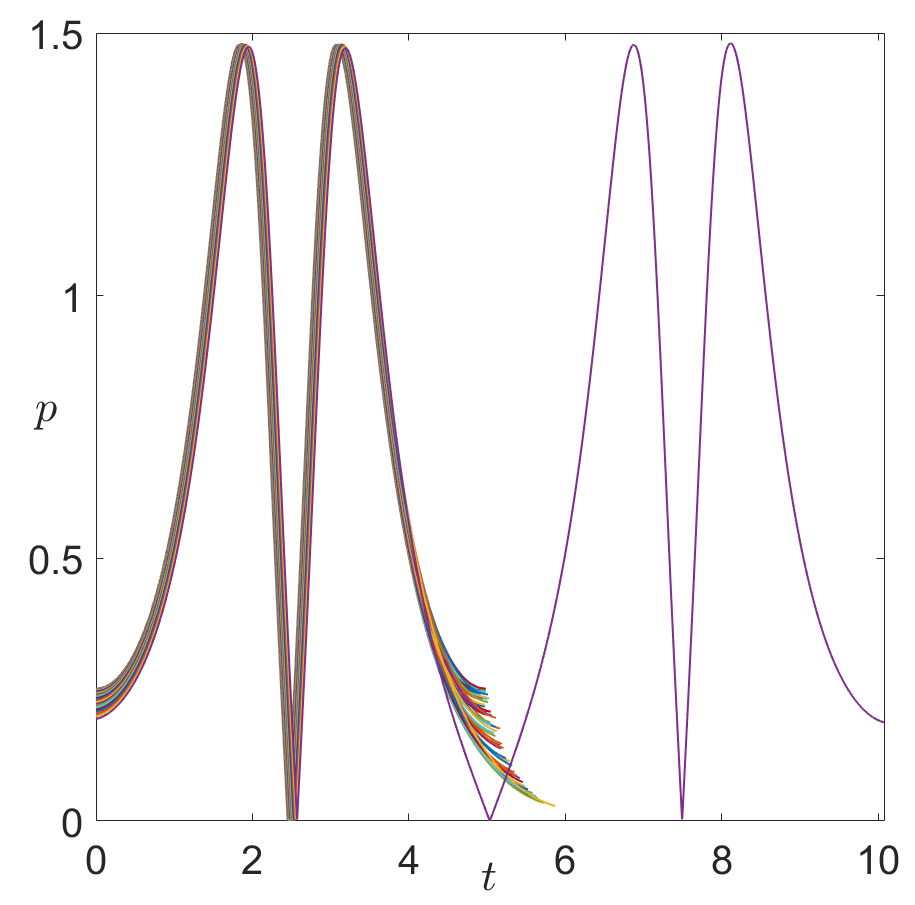} \\
		D)\includegraphics[scale = 0.2]{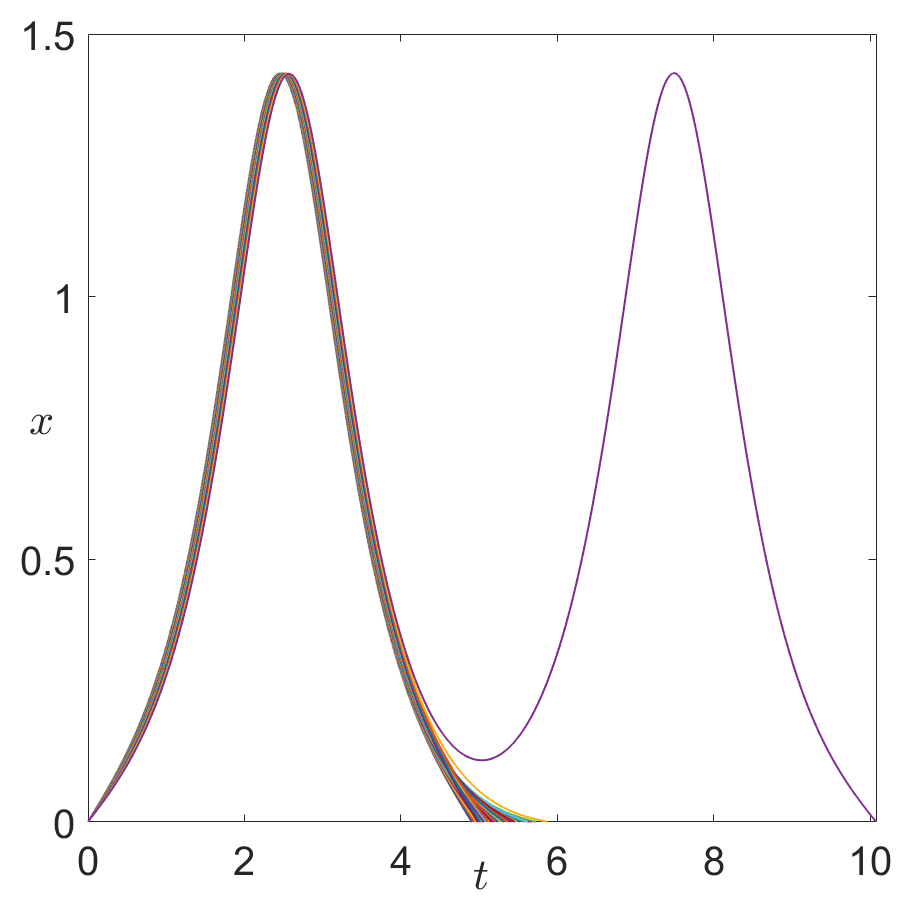}
		E)\includegraphics[scale = 0.2]{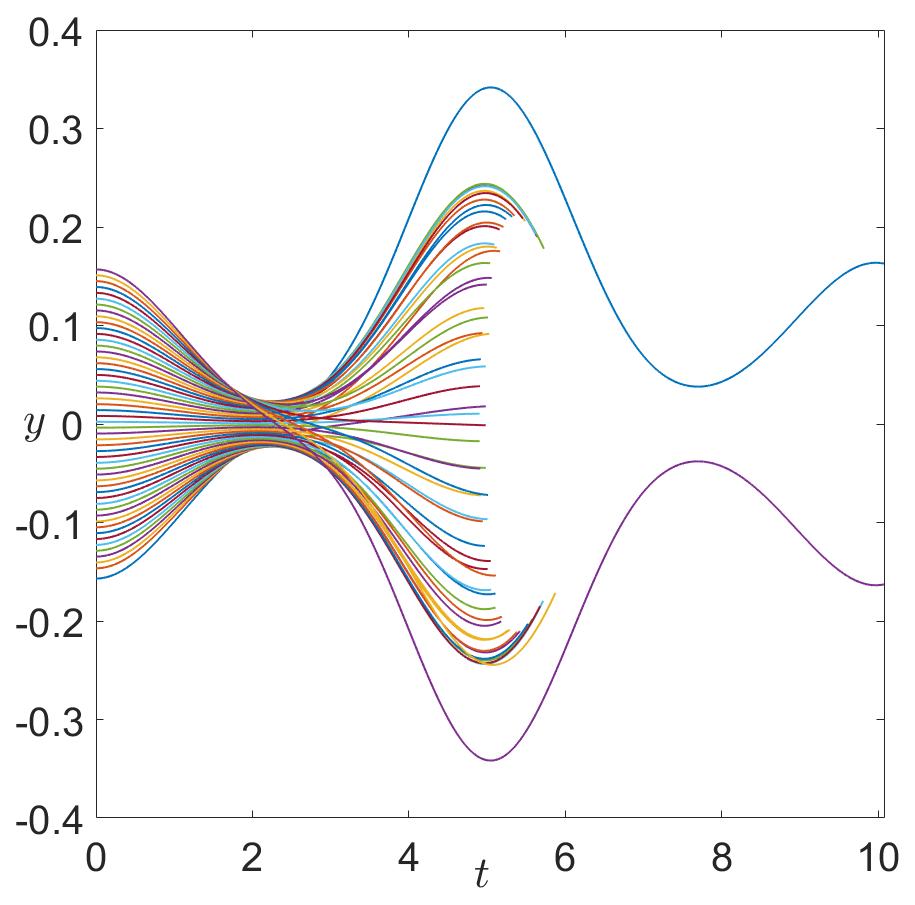}
	\end{center}
	\caption{Time evolution of the different components of a trajectory, displayed for all the recrossing trajectories obtained from an ensemble initialized at the configuration space line given by Eq. \eqref{eq:cs_slice}. The energy of the system is set to $\mathcal{H}_0 = 0.03$ and the location of the VRI point is at $x_i = 0.3265$. A) $p_x$ component of momentum; B) $p_y$ component of momentum; C) Total momentum; D) $x$ coordinate of the trajectory; E) $y$ coordinate of the trajectory.}
	\label{ensemble_timeEvol}
\end{figure}

\begin{figure}[htbp]
	\begin{center}
		\includegraphics[scale = 0.68]{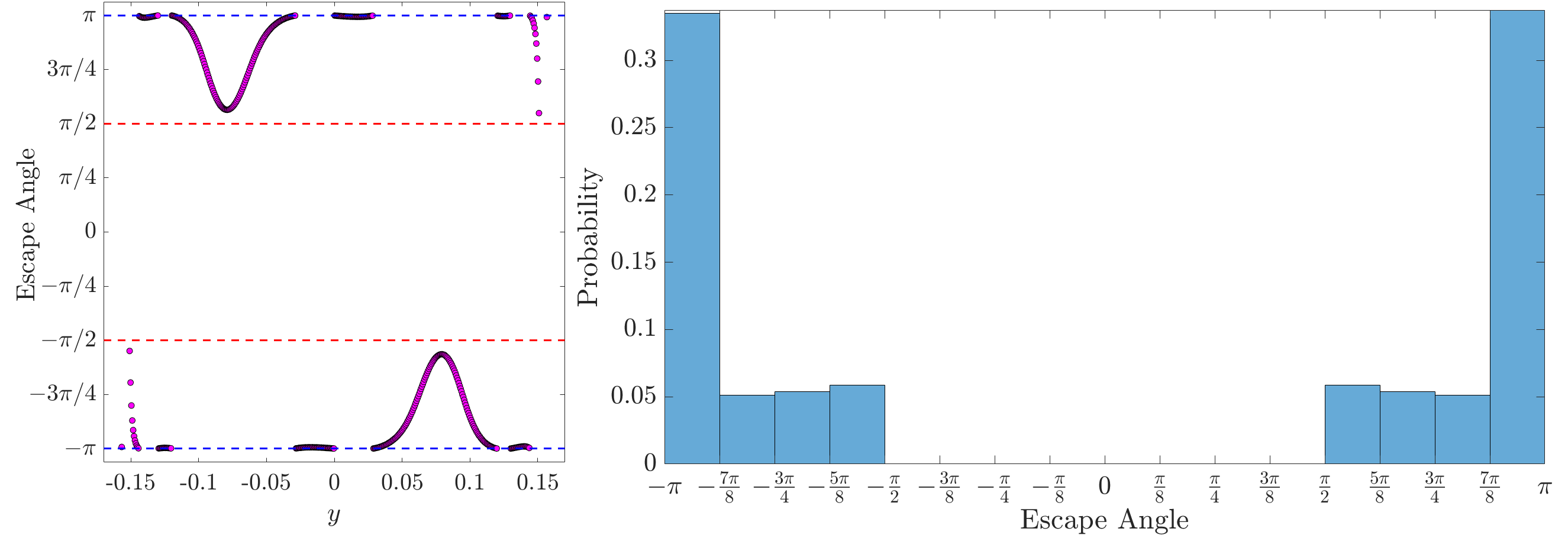}
	\end{center}
	\caption{Statistical distribution of the angles at which trajectories initialized on the configuration space line defined in Eq. \eqref{eq:cs_slice} recross the high energy saddle. We discard all those trajectories that enter either potential well on the right hand side of the PES before recrossing (escaping). The energy of the system is $\mathcal{H}_0 = 0.03$ and the location of the VRI point is at $x_i = 0.3265$. The $y$ coordinate represents the location of the initial condition along the configuration space line in Eq. \eqref{eq:cs_slice}.}
	\label{angle_ensemble_stats}
\end{figure}

From these simulations, it is also important to highlight that recrossing trajectories appear to have a tendency to preserve, to a certain extent, the 'directionality' at which they were initialized. What we mean by directionality is that the angle at which the trajectory is initialized and the angle at recrossing approximately  differ by $180^{\circ}$ (see the results displayed in Fig. \ref{angle_ensemble_stats} for a particular ensemble). This behavior of the ensemble of recrossing trajectories can be viewed and interpreted as some type of 'dynamical matching' mechanism \cite{carpenter1995,carpenter1998dynamic}. Interestingly, this effect on trajectories has also been recognized as relevant for chemical systems with PTSBs, such as in the Diels-Alder reaction \cite{singleton2009}. In order to address this question further, we consider the case where the energy is $\mathcal{H}_0 = 0.03$ and put the VRI point at the location $x_i = 0.3265$ for which the recrossing fraction gets its maximum value. In Fig. \ref{configSp_evol2} we depict all the recrossing trajectories projected onto configuration space, and we stop their evolution when they cross the $y$ axis. We have also overlaid on the PES, the value of the components of the force at every point of the configuration space. Two particular trajectories stand out in their evolution from the rest of the ensemble, and we have marked them in blue and magenta. These trajectories correspond to initial conditions whose trajectories get the 'closest' to entering the potential well regions. This is so, because they get very close to crossing the periodic trajectories that control access to the well region of the PES, and whose existence was demonstrated in previous work \cite{Agaoglou2020,katsanikas2020PRE} for this type of symmetric PESs with VRI points. If we analyze the time evolution of each of the components of the recrossing trajectories, see Fig. \ref{ensemble_timeEvol} we can see that the ensemble evolves in a coherent way, similar to the propagation of a wavepacket or a soliton-type solution. The 'limiting' trajectories that take longer to recross are the blue and magenta trajectories displayed in Fig. \ref{configSp_evol2}, which provide a 'natural' boundary to distinguish the trajectories that enter the well regions of the PES from those that recross. Notice also that the time evolution of their components also acts as an envelope for the time evolution of the components of all recrossing trajectories.

\begin{figure}[htbp]
	\begin{center}
		\includegraphics[scale = 0.68]{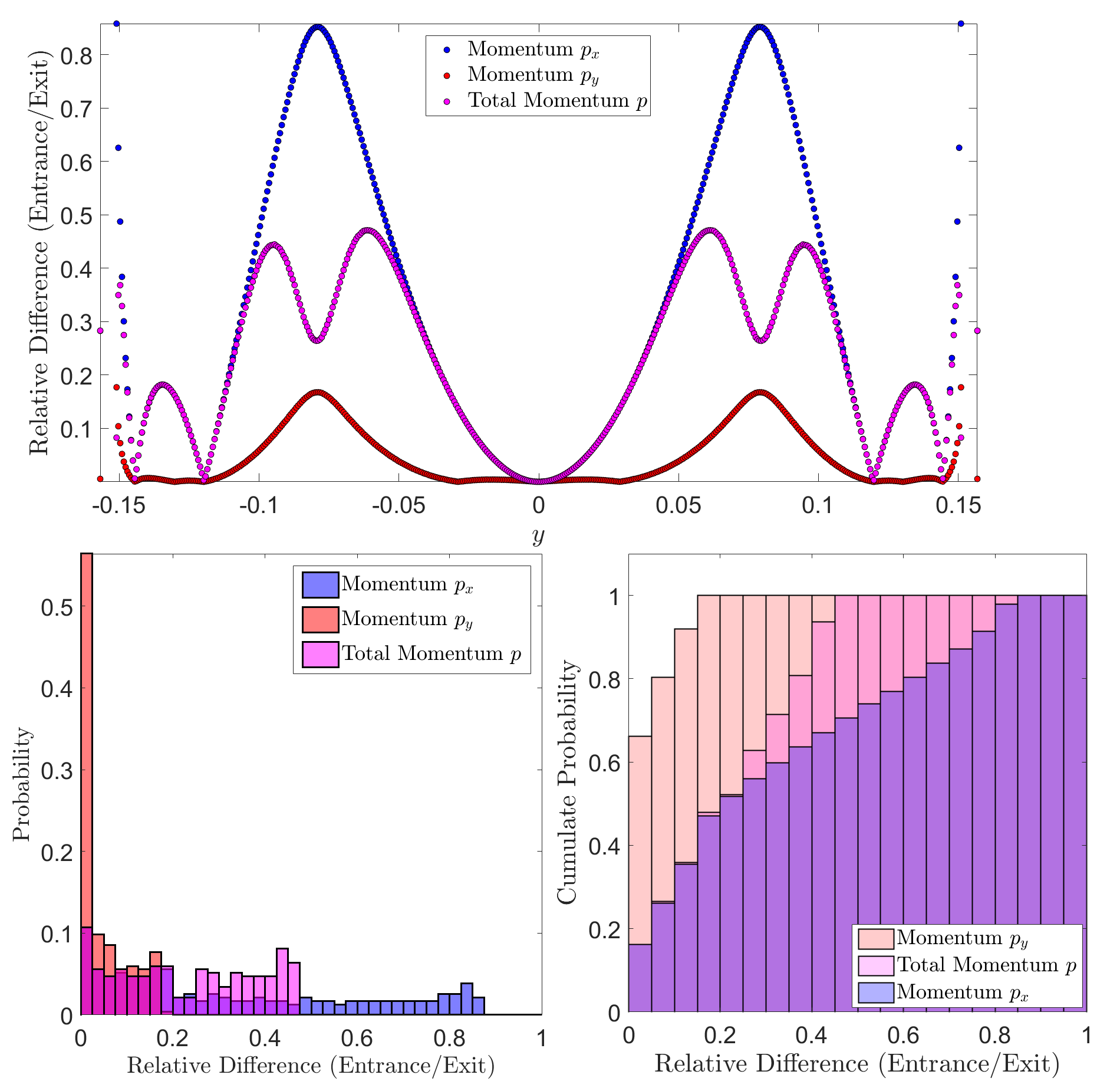}
	\end{center}
	\caption{Statistics for the relative difference in momentum values for recrossing trajectories initialized on the configuration space line defined in Eq. \eqref{eq:cs_slice}. The energy of the system is taken as $\mathcal{H}_0 = 0.03$ and the location of the VRI point is at $x_i = 0.3265$. The $y$ coordinate represents the location of the initial condition along the configuration space line in Eq. \eqref{eq:cs_slice}. The relative difference is calculated by subtracting the initial and final momentum values and dividing the result by the initial momentum. For the momentum $p_y$ we show the absolute difference, since the initial momentum $p_y = 0$ for all the trajectories.}
	\label{Mom_ensemble_stats}
\end{figure}

By performing a statistical analysis on the components of the recrossing trajectories at the instant in which they cross the high energy saddle at the origin, we can provide further evidence that they are exhibiting dynamical matching behavior. Recall that all initial conditions start from the line in Eq. \eqref{eq:cs_slice}, so that their initial momentum along the $y$ coordinate is zero, which means that they initially move horizontally. If we look at the momentum of the trajectories, see Fig. \ref{Mom_ensemble_stats}, in particular if we consider the relative difference between initial and final momentum values, we can see that many trajectories of the ensemble recross with a small value of $p_y$, and this is a clear indication that the directionality is preserved to a certain degree. We further check the angle at which the trajectories recross the $y$-axis (measured from the horizontal axis) to compare the horizontal deviation of the outgoing trajectory. Results are displayed in Fig. \ref{angle_ensemble_stats}, validating that directionality remains nearly horizontal at the point where they recross for many trajectories of the ensemble.  

\begin{figure}[htbp]
	\begin{center}
		A)\includegraphics[scale = 0.2]{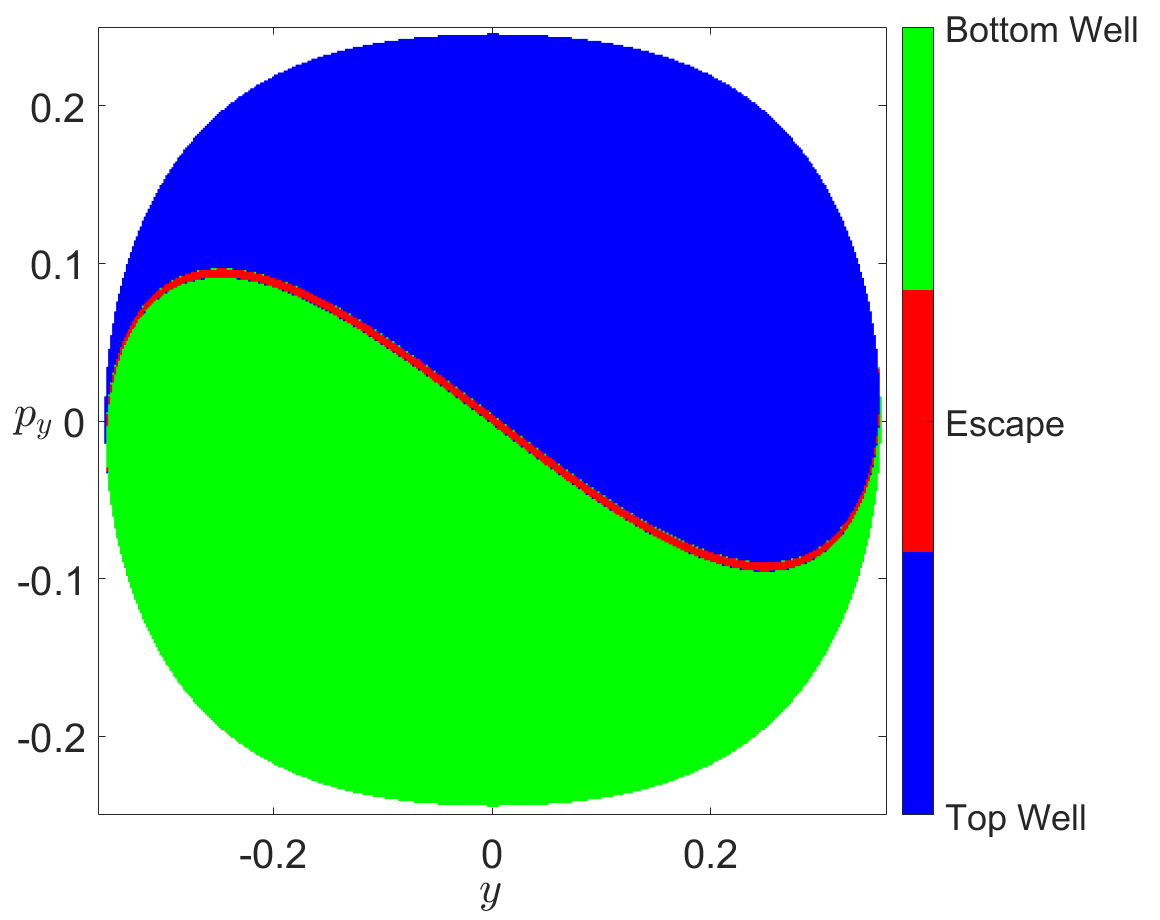}
		B)\includegraphics[scale = 0.2]{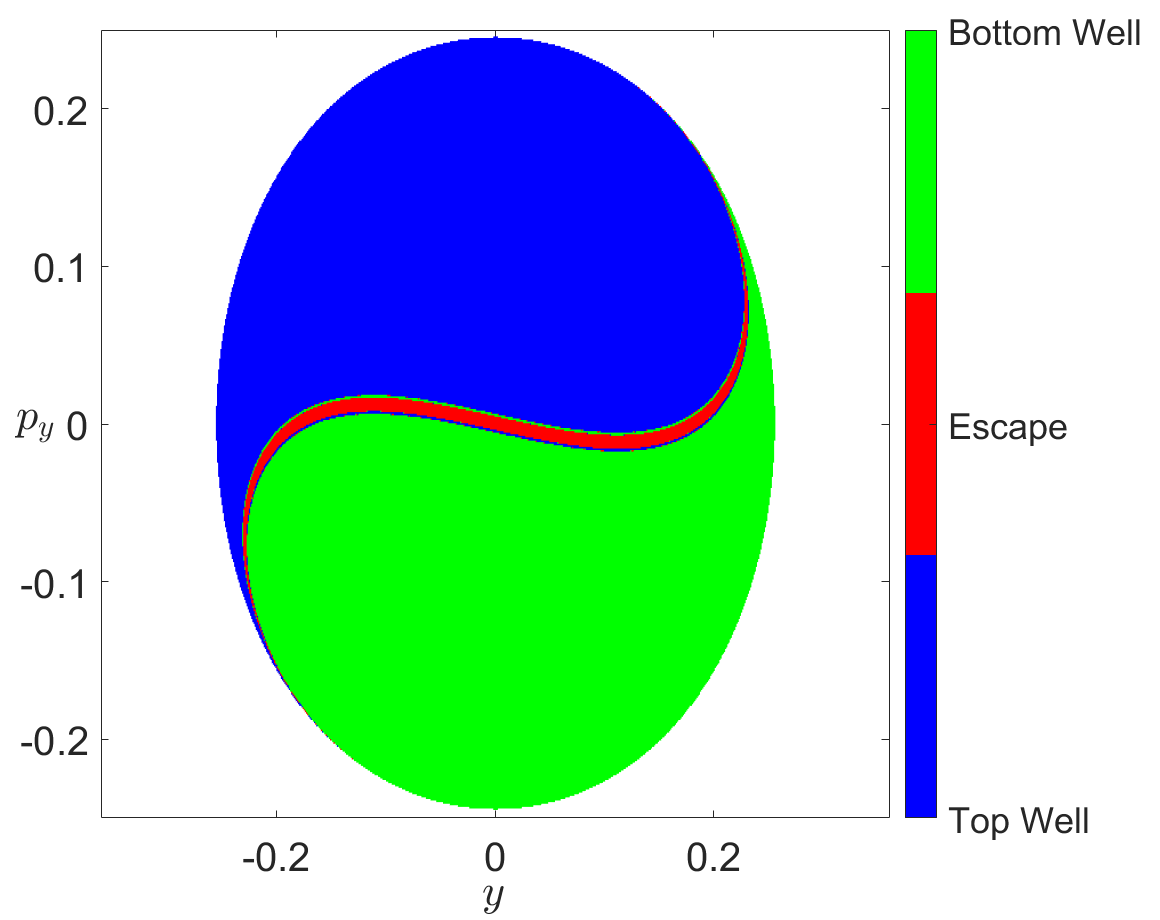} \\
		C)\includegraphics[scale = 0.2]{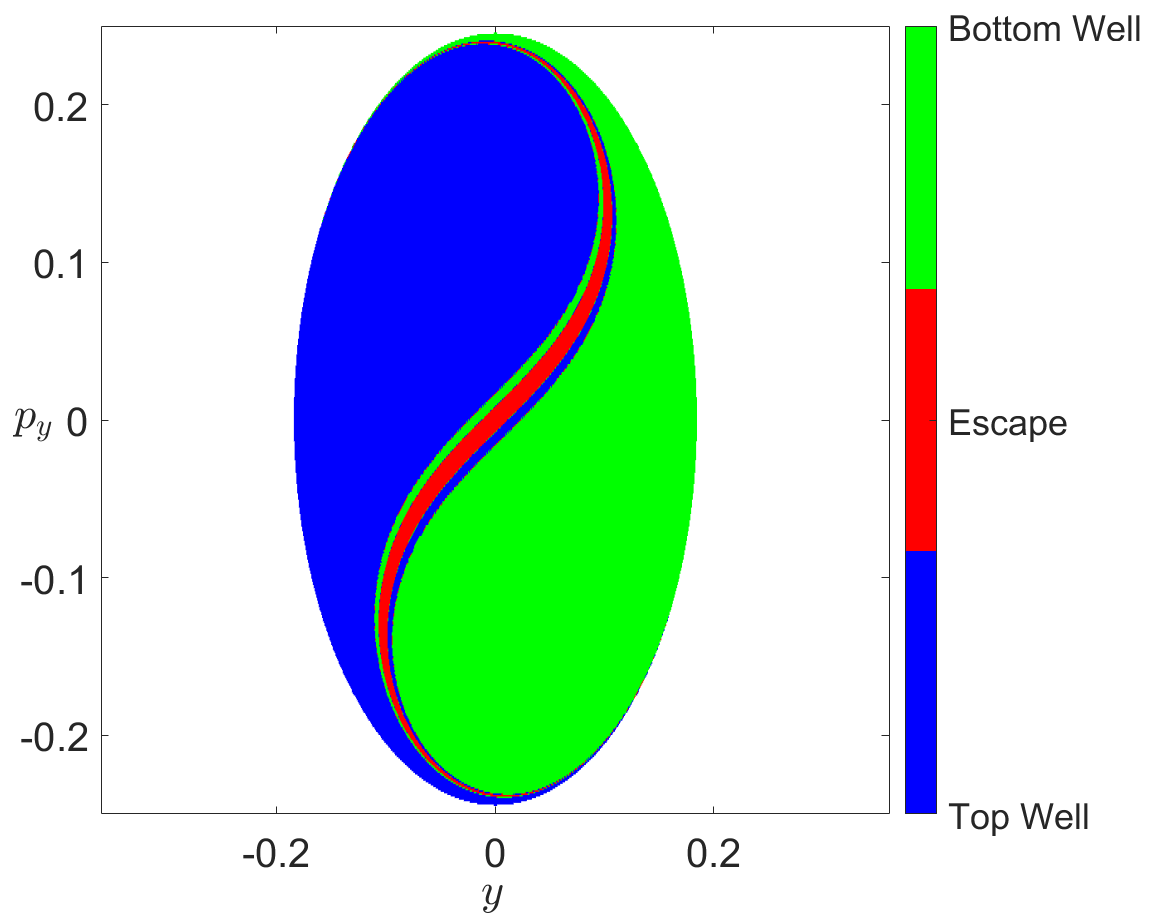}
		D)\includegraphics[scale = 0.2]{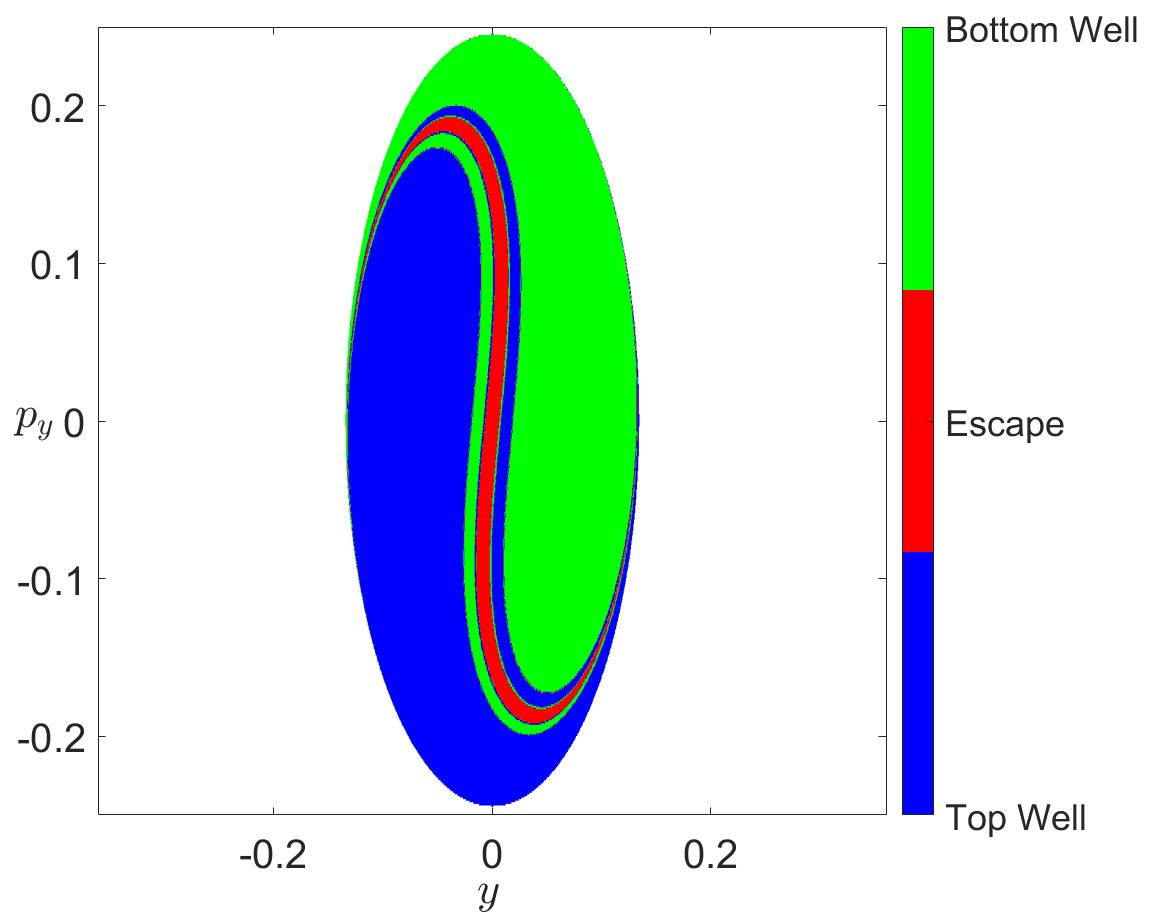}
	\end{center}
	\caption{Fate maps calculated on the phase space section described by Eq. \eqref{eq:ps_slice} for an energy of the system $\mathcal{H}_0 = 0.03$. Panels correspond to different locations of the VRI point. A) The VRI point is at $x_i = 0.1$; B) $x_i = 0.3$; C) $x_i = 0.5$; D) $x_i = 0.7$.}
	\label{fateMaps}
\end{figure}

We finish this work by studying the impact that the location of the VRI point has on the phase space structure of the Hamiltonian system. In order to address this question, we carry out the following simulation. First, we set the energy of the system to $\mathcal{H}_0 = 0.03$ and consider the two-dimensional phase space slice $\mathcal{P}(\mathcal{H}_0)$ described in Eq. \eqref{eq:ps_slice}. This plane is taken at the bottleneck region of the index-1 saddle equilibrium point at the origin, and we define a uniform grid of initial conditions on it. All those that satisfy the energy constraint are integrated until they recross the phase space plane, or enter the top/ bottom well regions (we stop them when they enter a circle of radius $R = 0.2$ centered about either potential well). Recall that the initial conditions on the slice $\mathcal{P}(\mathcal{H}_0)$ can have non-zero $p_y$ momentum. This means that  trajectories can start moving at an angle, which allows us to extend the analysis we already performed for the set of initial conditions in $\mathcal{C}(\mathcal{H}_0)$. We color-code the fate of the trajectories in the ensemble and produce what is known as a fate map, which is depicted in Fig. \ref{fateMaps} for different values of the VRI location. We can see from the plots that, as the VRI point gets farther away from the saddle at the origin, the area of the energetically feasible region of initial conditions decreases. Therefore, the VRI point affects the geometry of the bottleneck region about the high energy saddle. But most importantly, the location of the VRI point has a clear and distinctive influence on the geometry of the regions that correspond to trajectories with distinct dynamical fates. As the VRI point approaches the lower energy saddle that sits between the two potential wells, the regions get distorted and twisted, rotating similarly to a corkscrew mechanism. Moreover, the regions split into thin bands that organize into interlaced layers, giving rise to a fractal-like pattern. It is also important to highlight that in Fig. \ref{fateMaps} C), the region that corresponds to recrossing trajectories gets 'aligned' with the horizontal axis of the plot, that is, with the $p_y = 0$ line. This will give rise to a large recrossing fraction for an ensemble of trajectories initialized along that line, which reproduces and explains the results we obtained before when analyzing ensembles taken on the set $\mathcal{C}(\mathcal{H}_0)$.

\begin{figure}[htbp]
	\begin{center}
		A)\includegraphics[scale = 0.26]{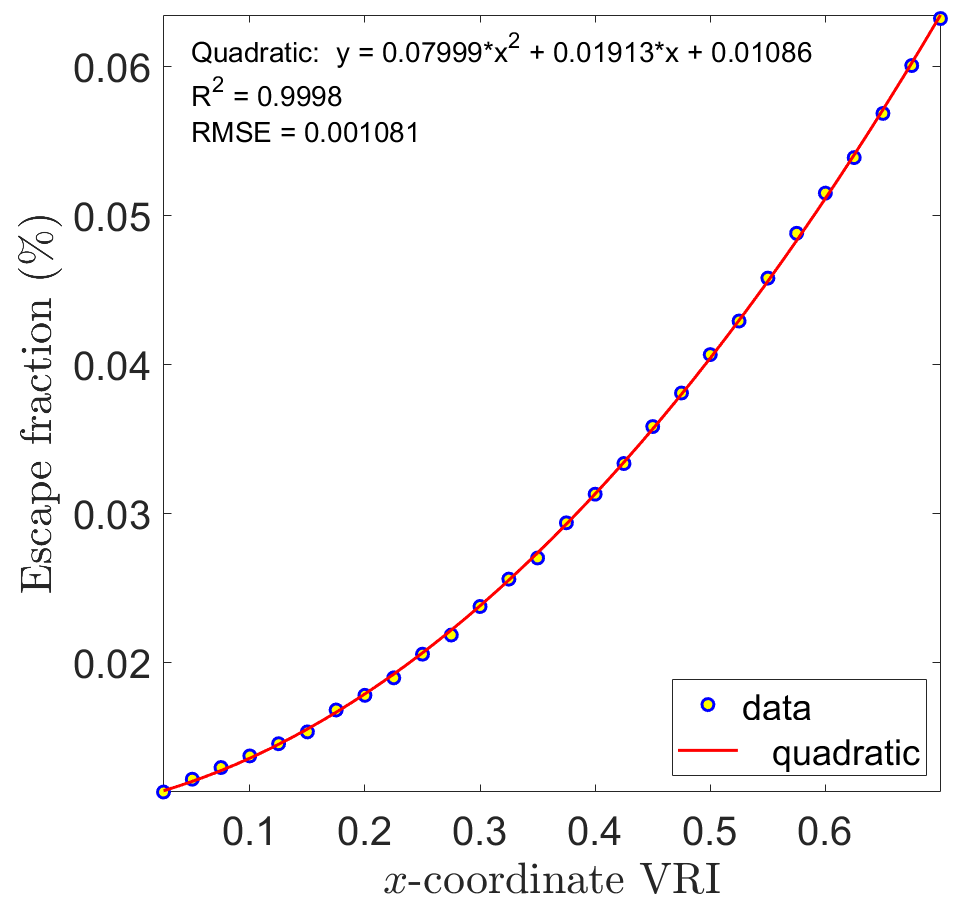}
		B)\includegraphics[scale = 0.26]{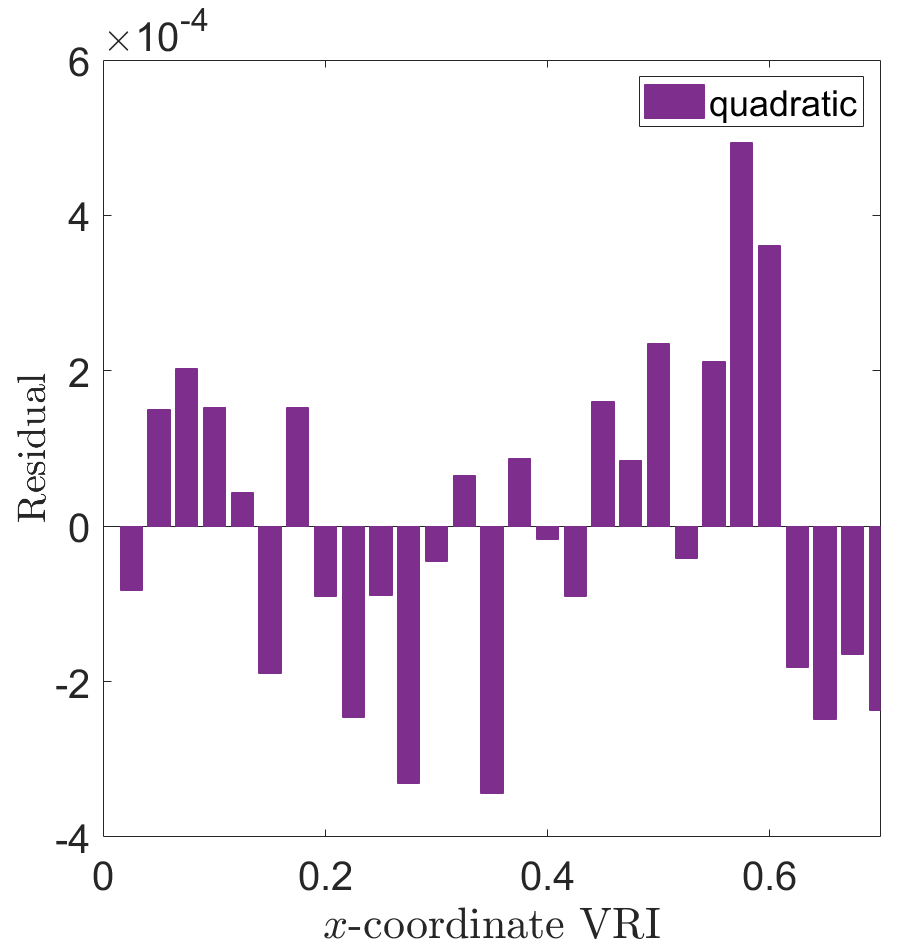} \\
		C)\includegraphics[scale = 0.26]{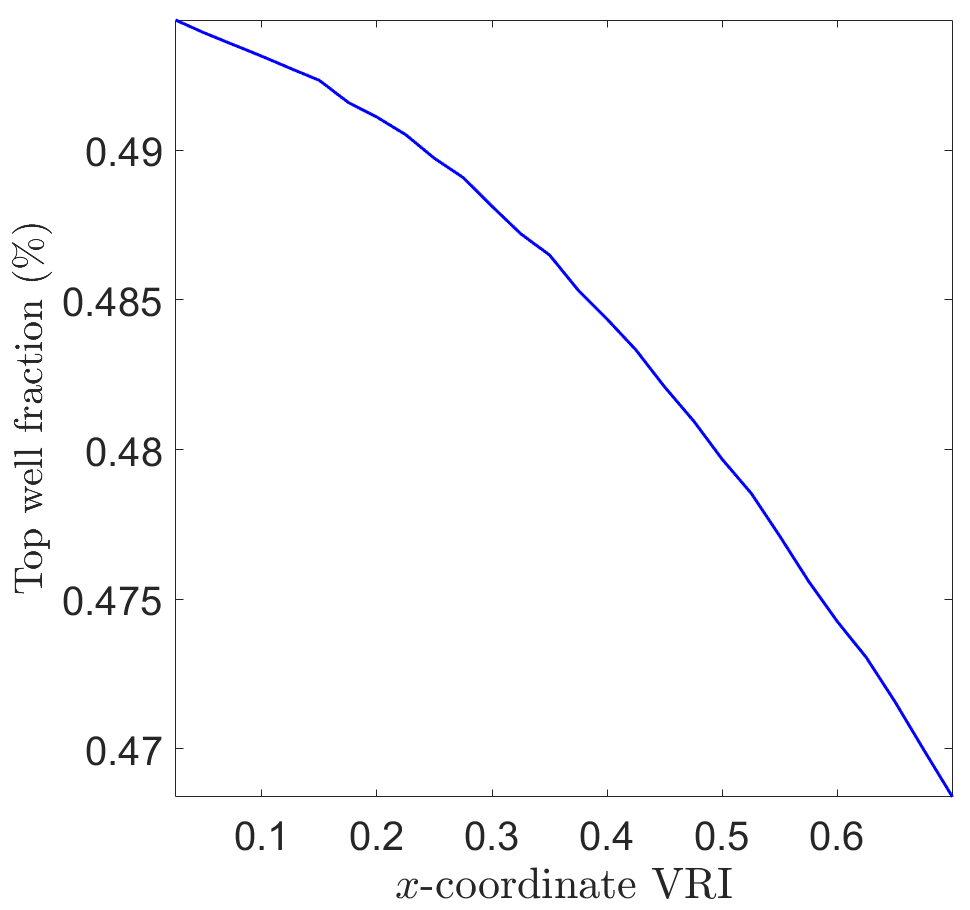}
		D)\includegraphics[scale = 0.26]{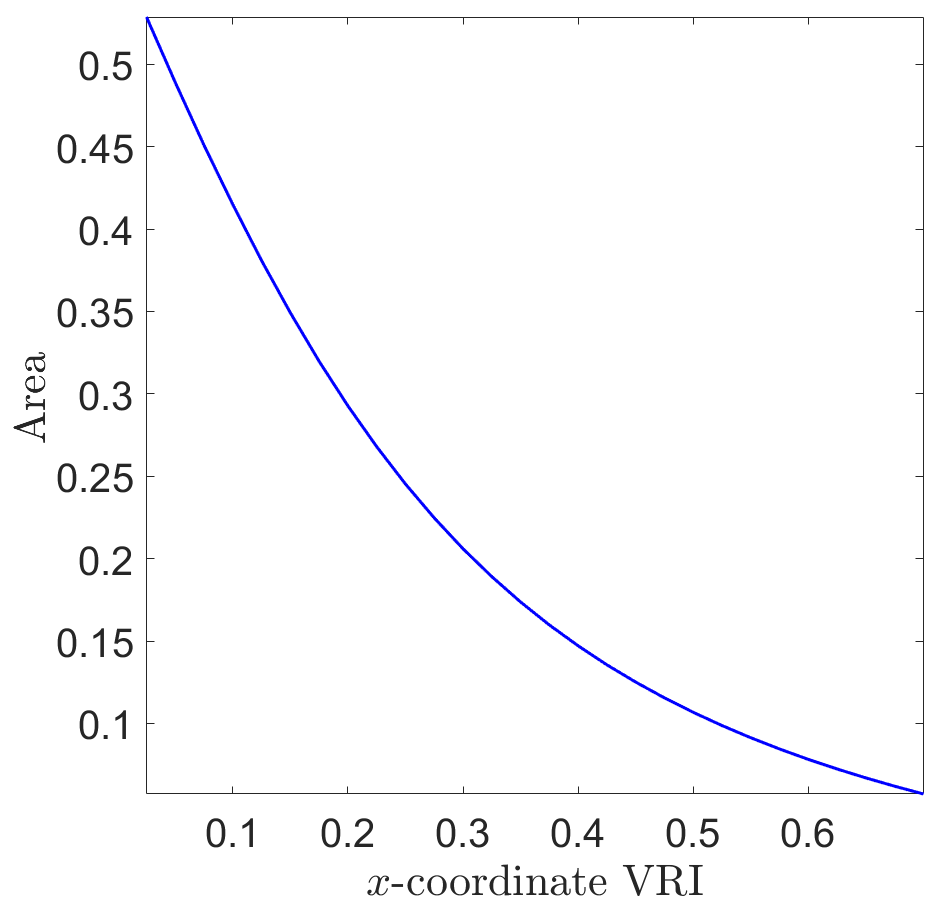}
	\end{center}
	\caption{Statistics for an ensemble of trajectories initialized on the phase space slice in Eq. \eqref{eq:ps_slice} with total energy $\mathcal{H}_0 = 0.03$. A) Fraction of recrossing trajectories as a function of the location of the VRI point. We have fitted the data obtained from numerical experiments to a quadratic model. B) Residuals of the quadratic fitting. C) Fraction of trajectories that enter the top well region of the PES. D) Variation of the area of the phase space slice $\mathcal{P}\left(\mathcal{H}_0\right)$ defined in Eq. \eqref{eq:ps_slice} with respect to the VRI location.}
	\label{escapeFraction}
\end{figure}

In order to quantify how the recrossing fraction of trajectories varies as a function of the VRI point location, we calculate the fate maps on the phase space section $\mathcal{P}(\mathcal{H}_0)$ for a range of values from $x_i = 0.025$ to $x_i = 0.7$ with a step of $\Delta x_i = 0.025$. The results of this simulation are presented in Fig. \ref{escapeFraction}. The fraction of recrossing trajectories is calculated by dividing the area of the slice by the area occupied by the region of recrossing trajectories (depicted in red in Fig. \ref{fateMaps}). We observe that, as the VRI point gets closer to the lower energy saddle of the PES, the recrossing fraction increases. Interestingly, the data obtained from the numerical experiments is nicely and accurately approximated by a quadratic law. 

\section{Conclusions}
\label{sec:conc}

In this work we have provided, for the first time, sound evidence that VRI points, despite not being equilibrium points of Hamilton's equations of motion, play an important role on the dynamics of trajectories, having a measurable impact on their dynamical fates. We have shown that for symmetric PESs, where the induced branching ratio is always $1:1$, the location of the VRI point controls the fraction of recrossing trajectories, that is, those trajectories that initially move across the high energy index-1 saddle and, after bouncing off the opposite wall of the PES, they go back to where they came from, without entering either of the wells separated by the lower energy saddle of the system. These trajectories do not give rise to the formation of products, and they reconvert to the original reactant configuration. Moreover, our numerical experiments point to the fact that recrossing trajectories have a tendency to behave as if they were experiencing some sort of ``dynamical matching'' mechanism, where the directionality along their evolution is preserved to a certain extent. 

Previous studies \cite{makrina2020cplett,katsanikas2020PRE,gg2020cplett} have highlighted the fundamental need for a phase space perspective to describe chemical reaction dynamics with PTSBs. In these works, it has been recognized that the dynamical mechanism in phase space, which determines selectivity in chemical reactions whose PES exhibits a PTSB region in their topography, is concerned with the existence of what is known as a heteroclinic connection between two unstable periodic orbits of the underlying Hamiltonian system. In this paper we have adopted a similar strategy, and our fate map analysis carried out in the phase space of the system has revealed two important features of the problem. One is that the location of the VRI point has a direct impact on the phase space structure, since the regions in the fate map corresponding to trajectories with distinct dynamical behavior get distorted and twisted, an effect similar to a corkscrew mechanism. On the other hand, given a fixed value for the energy of the system, the fraction of recrossing trajectories increases as the VRI point gets closer to the location of the lower energy saddle. Furthermore, there exists a quadratic relationship between the fraction of recrossing trajectories and the VRI point location. All the results obtained in the current work regarding the influence of the location of the VRI on the dynamical fate of trajectories point in the direction that the fraction of recrossing trajectories is controlled by a homoclinic connection. This geometrical structure is formed by the interaction between the stable and unstable manifolds of the unstable periodic orbit associated to the high energy saddle point at the origin. This phase space transport mechanism characterizes those trajectories that will recross the high energy saddle along their evolution. We will address these questions in further detail and explore this mechanism and its dynamical implications for chemical reactions in future work.

\section*{Acknowledgments}
The authors acknowledge the support of EPSRC Grant No. EP/P021123/1 and Office of Naval Research Grant No. N00014-01-1-0769.

\bibliography{vrilit}

\end{document}